\documentclass[final]{siamltex}%
\usepackage{eurosym}
\usepackage{algorithm}
\usepackage{algpseudocode}
\usepackage{amsfonts}
\usepackage{amsmath}
\usepackage{amssymb}
\usepackage{enumerate}
\usepackage{graphicx}
\usepackage{multirow}
\usepackage{url}
\usepackage{verbatim}
\usepackage[pdftex]{thumbpdf}
\usepackage[pdftex,
pdfstartview=FitH,bookmarks=true,bookmarksnumbered=true,bookmarksopen=true,hypertexnames=false,breaklinks=true,
colorlinks=true,  linkcolor=blue,anchorcolor=blue,
citecolor=blue,filecolor=blue,  menucolor=blue]{hyperref}
\usepackage{comment}%
\providecommand{\U}[1]{\protect\rule{.1in}{.1in}}

\providecommand{\U}[1]{\protect\rule{.1in}{.1in}}

\newtheorem{remark}[theorem]{Remark}
\newcommand*{\RE}
{_\text{Re}}
\newcommand*{\IM}
{_\text{Im}}
\begin{document}

\title{Stochastic Galerkin methods for linear stability analysis of systems with
parametric uncertainty\thanks{This work was supported by the U.\thinspace
\ S.\thinspace\ National Science Foundation under grant DMS1913201.}}
\author{Bed\v{r}ich Soused\'{\i}k\thanks{Department of Mathematics and Statistics,
University of Maryland, Baltimore County, 1000 Hilltop Circle, Baltimore,
MD~21250 (\texttt{sousedik@umbc.edu}) }
\and Kookjin Lee\thanks{ The School of Computing, Informatics, and Decision Systems
Engineering, Arizona State University, Tempe, AZ~85281
(\texttt{klee263@asu.edu}). } }
\maketitle

\begin{abstract}
We present a method for linear stability analysis of systems with parametric
uncertainty formulated in the stochastic Galerkin framework. Specifically, we
assume that for a model partial differential equation, the parameter is given
in the form of generalized polynomial chaos expansion. The stability analysis
leads to the solution of a stochastic eigenvalue problem, and we wish to
characterize the rightmost eigenvalue. We focus, in particular, on problems
with nonsymmetric matrix operators, for which the eigenvalue of interest may
be a complex conjugate pair, and we develop methods for their efficient
solution.
\textcolor{black}{These methods are based on inexact, line-search Newton iteration, which entails use of preconditioned GMRES.}
The method is applied to linear stability analysis of Navier\textendash Stokes
equation with stochastic viscosity, its accuracy is compared to that of Monte
Carlo and stochastic collocation, and the efficiency is illustrated by
numerical experiments.

\end{abstract}

\begin{keywords} linear stability, eigenvalue analysis, uncertainty quantification, spectral stochastic finite element method, Navier\textendash Stokes equation, preconditioning, stochastic Galerkin method \end{keywords}

\begin{AMS} 35R60, 60H15, 65F15, 65L07, 65N22, 65N25 \end{AMS}

\pagestyle{myheadings} \thispagestyle{plain} \markboth{B.\ SOUSED\'{I}K AND K. LEE}{STOCHASTIC GALERKIN METHODS FOR STABILITY ANALYSIS}

\section{Introduction}

\label{sec:introduction}The identification of instability in large-scale
dynamical system is important in a number of applications such as fluid
dynamics, epidemic models, pharmacokinetics, analysis of power systems and
power grid, or quantum mechanics and plasma physics.
\textcolor{black}{A steady solution~$u$ is stable,
if when in a transient simulation it is introduced with a small perturbation as initial data and the simulation reverts to~$u$,
and it is unstable otherwise. This is of fundamental importance since unstable solutions may lead to inexplicable dynamic behavior.
Linear stability analysis entails computing the rightmost eigenvalue of the Jacobian evaluated at~$u$,
and thus it leads to solution of, in general, large sparse generalized
eigenvalues problems see, e.g.~\cite{Cliffe-Spence-Taverner,Elman-2012-LII,Elman-2014-FEF,Govaerts,Kerner-1989-LSC,Schmid-Henningson} and the references therein.}
Typically, a complex pair of rightmost eigenvalues leads to a Hopf
bifurcation, and a real rightmost eigenvalue may lead to a pitchfork
bifurcation. The analysis is further complicated if the parameters in the
systems are functions of one or more random variables. This is quite common in
many real-world applications, since the precise values of model coefficients
or boundary conditions are often not known. A popular method for this type of
problems is Monte Carlo, which is known for its robustness but also slow
convergence. In this study, we use spectral stochastic finite element
methods~\cite{Ghanem-1991-SFE,LeMaitre-2010-SMU,Xiu-2010-NMS,Xiu-2002-WAP}
with the main focus on the so-called stochastic Galerkin method, for the
linear stability analysis of Navier\textendash Stokes equation with stochastic
viscosity. Specifically, we consider the parameterized\ viscosity given in the
form of generalized polynomial chaos (gPC) expansion. In the first step, we
apply the algorithms developed in~\cite{Lee-2019-LRS,Sousedik-2016-SGM}, see
also~\cite{Powell-2012-PSS}, to find a gPC\ expansion of the solution of the
Navier\textendash Stokes equation.
In the second step, we use the expansions of the solution and viscosity to set
up a generalized eigenvalue problem with a nonsymmetric matrix operator, and
in the assessment of linear stability of this problem we identify the
gPC\ expansions of the rightmost eigenvalue.
\textcolor{black}{The main contribution in this study is development of stochastic Galerkin method for nonsymmetric eigenvalue problems.
Our approach is based on inexact Newton iteration: the linear systems with Jacobian matrices are solved using GMRES,
for which we also develop several preconditioners.
The preconditioners are motivated by our prior work on (truncated) hierarchical preconditioning~\cite{Sousedik-2014-THP,Lee-2018-IMS},
see also~\cite{Bespalov-2021-TPS}.} For an overview of literature on solving
eigenvalue problems in the context of spectral stochastic finite element
methods we refer to~\cite{Andreev-2012-STA,Lee-2018-IMS,Sousedik-2016-ISI} and
the references therein.\ Recently, Hakula and Laaksonen~\cite{Hakula-2019-MSE}
studied crossing of eigenmodes in the stochastic parameter space,
and Elman and Su~\cite{Elman-2018-LRS} developed a low-rank inverse subspace
iteration.
\textcolor{black}{However, to the best of our knowledge, there are only a few
references addressing nonsymmetric stochastic eigenvalue problems: by Sarrouy
et al.~\cite{Sarrouy-2012-SAE,Sarrouy-2013-PPC}, but there is no discussion of
efficient solution strategies, and also by Sonday et
al.~\cite{Sonday-2011-EJG}, who studied distribution of eigenvalues for the
Jacobian in the context of stochastic Galerkin method.}
\textcolor{black}{Most recently, the authors with collaborators also compared surrogate learning strategies based on a sampling method, 
Gaussian process regression and a neural network in~\cite{Sousedik-2022-SLL}.}

A study of linear stability of Navier\textendash Stokes equation under
uncertainty was conducted by Elman and Silvester~\cite{Elman-2018-CME}. The
study was based on a judiciously chosen perturbation of the state variable and
a stochastic collocation method was used to characterize the rightmost
eigenvalue. Our approach here is different. We consider parametric uncertainty
(of the viscosity), and the solution strategy is based on the stochastic
Galerkin method. In fact, also the variant of the collocation method used here
is based on the stochastic Galerkin projection (sometimes called a
\emph{nonintrusive stochastic Galerkin method} in the literature,
see~\cite[Chapter~7]{Xiu-2010-NMS} for a discussion). From this perspective,
our study can be \textcolor{black}{viewed} as an extension of the setup
from~\cite{Elman-2012-LII} to problems with viscosity given in the form of
stochastic expansion and their efficient solution using stochastic Galerkin
method. However, more importantly, we illustrate that the inexact methods for
stochastic eigenvalue problems proposed recently by Lee and
Soused\'{\i}k~\cite{Lee-2018-IMS} can be also applied to problems with
nonsymmetric matrix operator\footnote{Specifically, the methods based on
inexact Newton iteration, since in our experience the stochastic Rayleigh
quotient and inverse iteration methods are less effective for nonsymmetric
problems.}. This in general allows to perform a linear stability analysis for
other types of problems as well.
\textcolor{black}{We do not address eigenvalue crossing here, which is a somewhat delicate task for gPC-based techniques.
We assume that the eigenvalue of interest is sufficiently separated from the rest of the spectrum and no crossing occurs. 
This is often the case for outliers and other eigenvalues that may be of interest in applications. 
A suitability of the algorithm we propose in this study can be assessed, e.g., by running first a low-fidelity (quasi-)Monte Carlo simulation.
We also note that from} our experience with the problem at hand indicates that
the rightmost eigenvalue remains relatively well separated from the rest of
the spectrum, and it is less prone to switching unlike the other eigenvalues, 
even with moderate values of coefficient of variation, which
in turn allows its efficient characterization by a gPC expansion.

The rest of the paper is organized as follows. In Section~\ref{sec:model} we
recall the basic elements of linear stability analysis and link it to an
eigenvalue problem for a specific model given by Navier\textendash Stokes
equation, in Section~\ref{sec:model-eig} we introduce the stochastic
eigenvalue problem, in Section~\ref{sec:sampling} we formulate the sampling
methods and in Section~\ref{sec:SG} the stochastic Galerkin method and inexact
Newton iteration for its solution, in Section~\ref{sec:NS} we apply the
algorithms to linear stability analysis of Navier\textendash Stokes equation
with stochastic viscosity, in Section~\ref{sec:numerical} we report the
results of numerical experiments, and finally in Section~\ref{sec:conclusion}
we summarize and conclude our work.

\section{Linear stability and deterministic model problem}

\label{sec:model}Following~\cite{Elman-2012-LII}, let us consider the
dynamical system
\begin{equation}
Mu_{t}=f(u,\nu), \label{eq:system-det}%
\end{equation}
where $f:%
\mathbb{R}
^{n}\times%
\mathbb{R}
{\color{black}\rightarrow}%
\mathbb{R}
^{n}$ is a nonlinear mapping, $u\in%
\mathbb{R}
^{n}$ is the state variable (velocity, pressure, temperature, deformation,
etc.), in the finite element setting $M\in%
\mathbb{R}
^{n\times n}$ is the mass matrix, and $\nu$ is a parameter. Let $u$ denote the
steady-state solution to~(\ref{eq:system-det}), i.e., $u_{t}=0$. We are
interested in the \emph{stability} of$~u$: if a small perturbation $\delta(0)$
is introduced to~$u$ at time $t=0$, does $\delta(t)$ grow with time, or does
it decay? For a fixed value of~$\nu$, linear stability of the steady-state
solution is determined by the spectrum of the eigenvalue problem
\begin{equation}
Jv=\lambda Mv, \label{eq:eig-det}%
\end{equation}
where $J=\frac{\partial f}{\partial u}(u(\nu),\nu)$ is the Jacobian matrix
of~$f$ evaluated at~$\nu$. The eigenvalues have a general form $\lambda
=\alpha+i\beta$, where $\alpha=\operatorname{Re}\lambda$ and $\beta
=\operatorname{Im}\lambda$. Then $e^{\lambda t}=e^{\alpha t}e^{i\beta t}$, and
since $\left\vert e^{\lambda t}\right\vert =\left\vert e^{\alpha t}\right\vert
$, there are \textcolor{black}{in general} two cases: if $\alpha<0$ the perturbation decays with time, and
if $\alpha>0$ the perturbation grows.
\textcolor{black}{We refer, e.g., to~\cite{Cliffe-Spence-Taverner,Govaerts} and the references therein for a detailed discussion.}
That is, if all the eigenvalues have strictly negative real part, then~$u$ is a
stable steady solution, and if some eigenvalues have nonnegative real part,
then$~u$ is unstable. Therefore, a change of stability can be detected by
monitoring the rightmost eigenvalues of~(\ref{eq:eig-det}). A steady-state
solution may lose its stability in one of two ways: either the rightmost
eigenvalue of~(\ref{eq:eig-det}) is real and passes through zero from negative
to positive as $\nu$ varies, or~(\ref{eq:eig-det}) has a complex pair of
rightmost eigenvalues and they cross the imaginary axis as $\nu$ varies, which
leads to a Hopf bifurcation with the consequent birth of periodic solutions
of~(\ref{eq:system-det}).

Consider a special case of~(\ref{eq:system-det}), the time-dependent
Navier\textendash Stokes equations governing viscous incompressible flow,
\begin{equation}
\begin{aligned} \vec{u}_{t}&=\nu \nabla ^{2}\vec{u}-\left( \vec{u}\cdot \nabla \right) \vec{u}-\nabla p, \\ 0&=\nabla \cdot \vec{u}, \end{aligned} \label{eq:NS-time}%
\end{equation}
subject to appropriate boundary and initial conditions in a bounded physical
domain$~D$, where $\nu$ is the kinematic viscosity, $\vec{u}$ is the velocity
and $p$ is the pressure. Properties of the flow are usually characterized by
the Reynolds number
\[
Re=\frac{UL}{\nu},
\]
where $U$ is the characteristic velocity and $L$ is the characteristic length.
For convenience, we will also sometimes refer to the Reynolds number instead
of the viscosity. Mixed finite element discretization of~(\ref{eq:NS-time})
gives the following Jacobian and the mass matrix, see~\cite{Elman-2012-LII}
and~\cite[Chapter~$8$]{Elman-2014-FEF} for more details,
\begin{equation}
J=\left[
\begin{array}
[c]{cc}%
F & B^{T}\\
B & 0
\end{array}
\right]  \in%
\mathbb{R}
^{n_{x}\times n_{x}},\qquad M=\left[
\begin{array}
[c]{cc}%
-G & 0\\
0 & 0
\end{array}
\right]  \in%
\mathbb{R}
^{n_{x}\times n_{x}}, \label{eq:matrices-det}%
\end{equation}
where $n_{x}=n_{u}+n_{p}$ is the number of velocity and pressure degrees of
freedom, respectively, $n_{u}>n_{p}$, $F\in%
\mathbb{R}
^{n_{u}\times n_{u}}$ is nonsymmetric, $B\in%
\mathbb{R}
^{n_{p}\times n_{u}}$ is the divergence matrix, and the velocity mass matrix
$G^{n_{u}\times n_{u}}$\ is symmetric positive definite. The matrices are
sparse and $n_{x}$ is typically large. We note that the mass matrix$~M$ is
singular, and (\ref{eq:eig-det}) has an infinite eigenvalue. As suggested
in~\cite{Cliffe-1994-EBM}, we replace the singular mass matrix$~M$ with the
nonsingular, shifted mass matrix
\begin{equation}
M_{\sigma}=\left[
\begin{array}
[c]{cc}%
-G & \sigma B^{T}\\
\sigma B & 0
\end{array}
\right]  , \label{eq:M-shifted}%
\end{equation}
which \textcolor{black}{is symmetric but in general indefinite, and it} maps
the infinite eigenvalues of~(\ref{eq:eig-det}) to $\sigma^{-1}$ and leaves the
finite ones unchanged. Then, the generalized eigenvalue
problem~(\ref{eq:eig-det}) can be transformed into an eigenvalue problem
\begin{equation}
Jv=\lambda M_{\sigma}v. \label{eq:eig-det-2}%
\end{equation}
The flow is considered stable if $\operatorname{Re}\lambda<0$, and we wish to
detect the onset of instability, that is to find when the rightmost eigenvalue
crosses the imaginary axis. Efficient methods for finding the rightmost pair
of complex eigenvalues of~(\ref{eq:eig-det}) (or~(\ref{eq:eig-det-2})) were
studied in~\cite{Elman-2012-LII}. Here, our goal is different. We consider
parametric uncertainty in the sense that the parameter $\nu\equiv\nu(\xi)$,
where $\xi$ is a set of random variables and which is given by the so-called
generalized polynomial chaos expansion. To this end, we first recast the
eigenvalue problem~(\ref{eq:eig-det-2}) in the spectral stochastic finite
element framework, then we show how to efficiently solve it, and finally we
apply the stability analysis to Navier\textendash Stokes equation with
stochastic viscosity.

\section{Stochastic eigenvalue problem}

\label{sec:model-eig}Let $\left(  \Omega,\mathcal{F},\mathcal{P}\right)  $ be
a complete probability space, that is $\Omega$ is the sample space with
$\sigma$-algebra$~\mathcal{F}$ and probability measure$~\mathcal{P}$. We
assume that the randomness in the mathematical model is induced by a vector
$\xi:\Omega\mapsto\Gamma\subset%
\mathbb{R}
^{m_{\xi}}$\ of independent 
random variables
$\xi_{1}(\omega),\dots,\xi_{m_{\xi}}(\omega)$, where $\omega\in\Omega$. Let
$\mathcal{B}(\Gamma)$ denote the Borel $\sigma$-algebra on$~\Gamma$ induced
by~$\xi$ and $\mu$ the induced measure. The expected value of the product of
measurable functions on$~\Gamma$ determines a Hilbert space $T_{\Gamma}\equiv
L^{2}\left(  \Gamma,\mathcal{B}(\Gamma),\mu\right)  $ with inner product
\begin{equation}
\left\langle u,v\right\rangle =\mathbb{E}\left[  uv\right]  =\int_{\Gamma
}u(\xi)v(\xi)d\mu(\xi), \label{eq:stoch-inner-prod}%
\end{equation}
where the symbol $\mathbb{E}$ denotes the mathematical expectation.

In computations we will use a finite-dimensional subspace $T_{p}\subset
T_{\Gamma}$ spanned by a set of multivariate polynomials $\left\{  \psi_{\ell
}(\xi)\right\}  $ that are orthonormal with respect to the density function
$\mu$, that is $\mathbb{E}\left[  \psi_{k}\psi_{\ell}\right]  =\delta_{k\ell}%
$, and $\psi_{1}$ is constant. This will be referred to as the gPC
basis~\cite{Xiu-2002-WAP}. The dimension of the space$~T_{p}$, depends on the
polynomial degree. For polynomials of total degree$~p$, the dimension is
$n_{\xi}=\binom{m_{\xi}+p}{p}$.

Suppose that we are given an affine expansion of a matrix operator, which may
correspond to the Jacobian matrix in~(\ref{eq:eig-det}), as%
\begin{equation}
K(\xi)=\sum_{\ell=1}^{n_{\nu}}K_{\ell}\psi_{\ell}(\xi), \label{eq:stoch-exp-A}%
\end{equation}
where each$~K_{\ell}\in\mathbb{R}^{n_{x}\times n_{x}}$ is a deterministic
matrix, and $K_{1}$ is the mean-value matrix $K_{1}=\mathbb{E}\left[
K(\xi)\right]  $. The representation~(\ref{eq:stoch-exp-A}) is obtained from
either Karhunen-Lo\`{e}ve or, more generally, a stochastic expansion of an
underlying random process; a specific example is provided in
Section~\ref{sec:numerical}.\ 

We are interested in a solution of the following stochastic eigenvalue
problem: find a specific eigenvalue~$\lambda(\xi)$ and possibly also the
corresponding eigenvector$~v(\xi)$, which satisfy in$~D$\ almost surely (a.s.)
the equation
\begin{equation}
K(\xi)v(\xi)=\lambda(\xi)M_{\sigma}v(\xi), \label{eq:param_eig}%
\end{equation}
where $K(\xi)\in\mathbb{R}^{n_{x}\times n_{x}}$ is a nonsymmetric matrix
operator, $M_{\sigma}\color{black}{\in\mathbb{R}^{n_{x}\times n_{x}}}$ is a deterministic mass matrix, $\lambda(\xi
)\in\mathbb{%
\mathbb{C}
}$ and $v(\xi)\in\mathbb{%
\mathbb{C}
}^{n_{x}}$ along with a normalization condition {\color{black}
\begin{equation}
\left(  v(\xi)\right)  ^{\ast}v(\xi)=\text{constant}, \label{eq:param_normal}%
\end{equation}
which is further specified in Section~\ref{sec:SG}.}

We will search for expansions of the eigenpair ${(\lambda(\xi),v(\xi))}$ in
the form
\begin{equation}
\lambda(\xi)=\sum_{k=1}^{n_{\xi}}\lambda_{k}\psi_{k}(\xi),\quad v(\xi
)=\sum_{k=1}^{n_{\xi}}v_{k}\psi_{k}(\xi), \label{eq:sol_mat}%
\end{equation}
where $\lambda_{k}\in\mathbb{%
\mathbb{C}
}$ and $v_{k}\in\mathbb{%
\mathbb{C}
}^{n_{x}}$ are coefficients corresponding to the basis$~\left\{  \psi
_{k}\right\}  $.
\textcolor{black}{We note that 
the series for $\lambda(\xi)$ in~(\ref{eq:sol_mat}) converges for $n_\xi \rightarrow \infty$ in the space~$T_{\Gamma}$ under the assumption that
the gPC polynomials provide its orthonormal basis and provided that~$\lambda(\xi)$ has finite second moments 
see, e.g.,~\cite{Ernst-2012-CGP} or~\cite[Chapter~5]{Xiu-2010-NMS}. 
Convergence analysis of this approximation 
for self-adjoint problems can be found in~\cite{Andreev-2012-STA}.}

\section{Monte Carlo and stochastic collocation}

\label{sec:sampling}Both Monte Carlo and stochastic collocation are based on
sampling. The coefficients are defined by a discrete projection
\begin{equation}
\lambda_{k}=\left\langle \lambda,\psi_{k}\right\rangle ,\quad v_{k}%
=\left\langle v,\psi_{k}\right\rangle ,\qquad k=1,\dots,n_{\xi}.
\label{eq:sol_mat-proj}%
\end{equation}
The evaluations of coefficients in~(\ref{eq:sol_mat-proj}) entail solving a
set of independent deterministic eigenvalue problems at a set of sample
points$~\left\{  \xi^{(q)}\right\}  $, $q=1,\dots,n_{MC}$ or $n_{q}$,
\begin{equation}
K(\xi^{(q)})v(\xi^{(q)})=\lambda\left(  \xi^{(q)}\right)  M_{\sigma}v\left(
\xi^{(q)}\right)  . \label{eq:param_eig_sampling}%
\end{equation}
In the Monte Carlo method, the sample points$~\xi^{(q)}$, $q=1,\dots,n_{MC}%
,$\ are generated randomly following the distribution of the random
variables$~\xi_{i}$, $i=1,\dots,m_{\xi}$, and moments of solution are computed
by ensemble averaging. In addition, the coefficients in~(\ref{eq:sol_mat})
\textcolor{black}{could} be computed
\textcolor{black}{using Monte Carlo integration} as\footnote{\textcolor{black}{In numerical
experiments, we avoid this approximation of the gPC coefficients and directly work with the sampled
quantities.}}%
\[
\lambda_{k}=\frac{1}{n_{MC}}\sum_{q=1}^{n_{MC}}\lambda(\xi^{(q)})\psi
_{k}\left(  \xi^{(q)}\right)  ,\qquad v_{mk}=\frac{1}{n_{MC}}\sum
_{q=1}^{n_{MC}}v(x_{m},\mathbb{\xi}^{(q)})\psi_{k}(\mathbb{\xi}^{(q)}).
\]
For stochastic collocation, which is used here as the so-called nonintrusive
stochastic Galerkin method, the sample points$~\xi^{(q)}$, $q=1,\dots,n_{q},$
consist of a predetermined set of \emph{collocation points}, and the
coefficients $\lambda_{k}$ and $v_{k}$\ in the expansions~(\ref{eq:sol_mat})
are determined by evaluating~(\ref{eq:sol_mat-proj}) in the sense
of~(\ref{eq:stoch-inner-prod}) using numerical quadrature as
\[
\lambda_{k}=\sum_{q=1}^{n_{q}}\lambda(\mathbb{\xi}^{(q)})\psi_{k}(\mathbb{\xi
}^{(q)})w^{(q)},\qquad v_{mk}=\sum_{q=1}^{n_{q}}v(x_{m},\mathbb{\xi}%
^{(q)})\psi_{k}(\mathbb{\xi}^{(q)})w^{(q)},
\]
where $\mathbb{\xi}^{(q)}$ are the quadrature (collocation) points and
$w^{(q)}$ are quadrature weights. Details of the rule we use in our numerical
experiments are discussed in Section~\ref{sec:numerical}, and we refer to
monograph~\cite{LeMaitre-2010-SMU}\ for a detailed discussion of quadrature rules.

\section{Stochastic Galerkin method and Newton iteration}

\label{sec:SG}The stochastic Galerkin method is based on the projection%
\begin{align}
\left\langle Kv,\psi_{k}\right\rangle  &  =\left\langle \lambda M_{\sigma
}v,\psi_{k}\right\rangle , & k  &  =1,\dots,n_{\xi},\label{eq:SG-eig}\\
\left\langle v^{{\color{black}\ast}}v,\psi_{k}\right\rangle  &  = {\color{black}\text{const} \cdot} \delta_{k1}, &
k  &  =1,\dots,n_{\xi}, \,\,\,  {\color{black}\text{const} \in \mathbb{R}}. \label{eq:SG-normal}%
\end{align}
Let us introduce the notation
\begin{equation}
\left[  H_{\ell}\right]  _{jk}=h_{\ell,jk},\quad h_{\ell,jk}\equiv
\mathbb{E}\left[  \psi_{\ell}\psi_{j}\psi_{k}\right]  ,\qquad\ell
=1,\dots,n_{\nu},\quad j,k=1,\dots,n_{\xi}. \label{eq:H}%
\end{equation}
In order to formulate an efficient algorithm for eigenvalue
problem~(\ref{eq:param_eig}) with nonsymmetric matrix operator using the
stochastic Galerkin formulation, we introduce a separated representation of
the eigenpair using real and imaginary parts,
\[
v(\xi)=v_{\operatorname{Re}}(\xi)+{\mathrm{i}\mkern1mu}v_{\operatorname{Im}%
}(\xi),\quad\lambda(\xi)=\lambda_{\operatorname{Re}}(\xi)+{\mathrm{i}%
\mkern1mu}\lambda_{\operatorname{Im}}(\xi),
\]
where $v_{\operatorname{Re}}(\xi),v_{\operatorname{Im}}(\xi)\in\mathbb{R}%
^{n_{x}}$ and $\lambda_{\operatorname{Re}}(\xi),\lambda_{\operatorname{Im}%
}(\xi)\in\mathbb{R}$. Then, replacing $v(\xi)$ and $\lambda(\xi)$ in
eigenvalue problem~(\ref{eq:param_eig}) results in
\begin{equation}
K(\xi)(v_{\operatorname{Re}}(\xi)+{\mathrm{i}\mkern1mu}v_{\operatorname{Im}%
}(\xi))=\left(  \lambda_{\operatorname{Re}}(\xi)+{\mathrm{i}\mkern1mu}%
\lambda_{\operatorname{Im}}(\xi)\right)  M_{\sigma}(v_{\operatorname{Re}}%
(\xi)+{\mathrm{i}\mkern1mu}v_{\operatorname{Im}}(\xi)).
\label{eq:sep_param_eig}%
\end{equation}
Expanding the terms in~(\ref{eq:sep_param_eig}) and collecting the real and
imaginary parts yields a system of equations that can be written 
\textcolor{black}{in a separated representation} as
\begin{equation}
\left[
\begin{array}
[c]{c}%
K(\xi)v_{\operatorname{Re}}(\xi)\\
K(\xi)v_{\operatorname{Im}}(\xi)
\end{array}
\right]  =%
\begin{bmatrix}
\lambda_{\operatorname{Re}}(\xi)M_{\sigma} & -\lambda_{\operatorname{Im}}%
(\xi)M_{\sigma}\\
\lambda_{\operatorname{Im}}(\xi)M_{\sigma} & \lambda_{\operatorname{Re}}%
(\xi)M_{\sigma}%
\end{bmatrix}
\left[
\begin{array}
[c]{c}%
v_{\operatorname{Re}}(\xi)\\
v_{\operatorname{Im}}(\xi)
\end{array}
\right]  , \label{eq:sep_param_eig2}%
\end{equation}
and the \textcolor{black}{normalization condition corresponding to the separated 
representation in~(\ref{eq:sep_param_eig2}) 
is} taken as
\begin{equation}
v_{\operatorname{Re}}(\xi)^{T}v_{\operatorname{Re}}(\xi)=1,\quad
v_{\operatorname{Im}}(\xi)^{T}v_{\operatorname{Im}}(\xi)=1.
\label{eq:sep_norm}%
\end{equation}
Now, we seek expansions of type~(\ref{eq:sol_mat}) for $v_{\operatorname{Re}%
}(\xi)$, $v_{\operatorname{Im}}(\xi)$, $\lambda_{\operatorname{Re}}(\xi)$, and
$\lambda_{\operatorname{Im}}(\xi)$, that is
\begin{equation}
v_{s}(\xi)=(\Psi(\xi)^{T}\otimes I_{n_{x}})\bar{v}_{s},\quad\lambda_{s}%
(\xi)=\Psi(\xi)^{T}\bar{\lambda}_{s},\quad s=\operatorname{Re}%
,\operatorname{Im}, \label{eq:param_sol}%
\end{equation}
where the symbol~$\otimes$ denotes the Kronecker product, $\Psi(\xi)=[\psi
_{1}(\xi),\ldots,\psi_{n_{\xi}}(\xi)]^{T}$, $\bar{\lambda}_{s}=[\lambda
_{s,1},\ldots,\lambda_{s,n_{\xi}}]^{T}\in\mathbb{R}^{n_{\xi}}$, and $\bar
{v}_{s}=[v_{s,1}^{T},\ldots,v_{s,n_{\xi}}^{T}]^{T}\in\mathbb{R}^{n_{x}n_{\xi}%
}$ for $s=\operatorname{Re},\operatorname{Im}$.

Let us consider expansions~(\ref{eq:param_sol}) as approximations to the
solution of~(\ref{eq:sep_param_eig2})--(\ref{eq:sep_norm}). Then, we can write
the residual of~(\ref{eq:sep_param_eig2}) as
\[
\widetilde{F}(\bar{v}_{\operatorname{Re}},\bar{v}_{\operatorname{Im}}%
,\bar{\lambda}_{\operatorname{Re}},\bar{\lambda}_{\operatorname{Im}})=%
\begin{bmatrix}
\left(  \Psi(\xi)^{T}\!\otimes\!K(\xi)\right)  \bar{v}_{\operatorname{Re}%
}-\!\left(  \bar{\lambda}_{\operatorname{Re}}^{T}\Psi(\xi)\Psi(\xi
)^{T}\!\otimes\!M_{\sigma}\right)  \bar{v}_{\operatorname{Re}}\!+\!\left(
\bar{\lambda}_{\operatorname{Im}}^{T}\Psi(\xi)\Psi(\xi)^{T}\!\otimes
\!M_{\sigma}\right)  \bar{v}_{\operatorname{Im}}\\
\vspace{-2.5mm}\\
\left(  \Psi(\xi)^{T}\!\otimes\!K(\xi)\right)  \bar{v}_{\operatorname{Im}%
}\!-\!\left(  \bar{\lambda}_{\operatorname{Im}}^{T}\Psi(\xi)\Psi(\xi
)^{T}\!\otimes\!M_{\sigma}\right)  \bar{v}_{\operatorname{Re}}\!-\!\left(
\bar{\lambda}_{\operatorname{Re}}^{T}\Psi(\xi)\Psi(\xi)^{T}\!\otimes
\!M_{\sigma}\right)  \bar{v}_{\operatorname{Im}}%
\end{bmatrix}
,
\]
and the residual of~(\ref{eq:sep_norm}) as
\[
\widetilde{G}(\bar{v}_{\operatorname{Re}},\bar{v}_{\operatorname{Im}})=%
\begin{bmatrix}
\bar{v}_{\operatorname{Re}}^{T}\left(  \Psi(\xi)\Psi(\xi)^{T}\otimes I_{n_{x}%
}\right)  \bar{v}_{\operatorname{Re}}-1\\
\vspace{-2.5mm}\\
\bar{v}_{\operatorname{Im}}^{T}\left(  \Psi(\xi)\Psi(\xi)^{T}\otimes I_{n_{x}%
}\right)  \bar{v}_{\operatorname{Im}}-1
\end{bmatrix}
.
\]
Imposing the stochastic Galerkin orthogonality conditions~(\ref{eq:SG-eig})
and~(\ref{eq:SG-normal}) to$~\widetilde{F}$ and$~\widetilde{G}$, respectively,
we get a system of nonlinear equations
\begin{equation}
r(\bar{v}_{\operatorname{Re}},\bar{v}_{\operatorname{Im}},\bar{\lambda
}_{\operatorname{Re}},\bar{\lambda}_{\operatorname{Im}})=%
\begin{bmatrix}
F(\bar{v}_{\operatorname{Re}},\bar{v}_{\operatorname{Im}},\bar{\lambda
}_{\operatorname{Re}},\bar{\lambda}_{\operatorname{Im}})\\
G(\bar{v}_{\operatorname{Re}},\bar{v}_{\operatorname{Im}})
\end{bmatrix}
={0}\in\mathbb{R}^{2(n_{x}+1)n_{\xi}}, \label{eq:nonlinear_res}%
\end{equation}
where
\[
F(\bar{v}_{\operatorname{Re}},\bar{v}_{\operatorname{Im}},\bar{\lambda
}_{\operatorname{Re}},\bar{\lambda}_{\operatorname{Im}})=%
\begin{bmatrix}
\mathbb{E}[\Psi\Psi^{T}\!\otimes\!K]\bar{v}_{\operatorname{Re}}\!-\!\mathbb{E}%
[(\bar{\lambda}_{\operatorname{Re}}{}^{T}\Psi)\Psi\Psi^{T}\!\otimes
\!M_{\sigma}]\bar{v}_{\operatorname{Re}}\!+\!\mathbb{E}[(\bar{\lambda
}_{\operatorname{Im}}{}^{T}\Psi)\Psi\Psi^{T}\!\otimes\!M_{\sigma}]\bar
{v}_{\operatorname{Im}}\\
\vspace{-2.5mm}\\
\mathbb{E}[\Psi\Psi^{T}\!\otimes\!K]\bar{v}_{\operatorname{Im}}\!-\!\mathbb{E}%
[(\bar{\lambda}_{\operatorname{Im}}{}^{T}\Psi)\Psi\Psi^{T}\!\otimes
\!M_{\sigma}]\bar{v}_{\operatorname{Re}}\!-\!\mathbb{E}[(\bar{\lambda
}_{\operatorname{Re}}{}^{T}\Psi)\Psi\Psi^{T}\!\otimes\!M_{\sigma}]\bar
{v}_{\operatorname{Im}}%
\end{bmatrix}
,
\]
and
\[
G(\bar{v}_{\operatorname{Re}},\bar{v}_{\operatorname{Im}})=%
\begin{bmatrix}
\mathbb{E}[\Psi\otimes((\bar{v}_{\operatorname{Re}}{}^{T}(\Psi\Psi^{T}\otimes
I_{n_{x}})\bar{v}_{\operatorname{Re}})-1)]\\
\vspace{-2.5mm}\\
\mathbb{E}[\Psi\otimes((\bar{v}_{\operatorname{Im}}{}^{T}(\Psi\Psi^{T}\otimes
I_{n_{x}})\bar{v}_{\operatorname{Im}})-1)]
\end{bmatrix}
.
\]

We will use Newton iteration to solve system of nonlinear
equations~(\ref{eq:nonlinear_res}). The Jacobian matrix$~DJ(\bar
{v}_{\operatorname{Re}},\bar{v}_{\operatorname{Im}},\bar{\lambda
}_{\operatorname{Re}},\bar{\lambda}_{\operatorname{Im}})$
of~(\ref{eq:nonlinear_res}) can be written as
\[
DJ(\bar{v}_{\operatorname{Re}},\bar{v}_{\operatorname{Im}},\bar{\lambda
}_{\operatorname{Re}},\bar{\lambda}_{\operatorname{Im}})=%
\begin{bmatrix}
\frac{\partial F}{\partial\bar{v}_{\operatorname{Re}}} & \frac{\partial
F}{\partial\bar{v}_{\operatorname{Im}}} & \frac{\partial F}{\partial
\bar{\lambda}_{\operatorname{Re}}} & \frac{\partial F}{\partial\bar{\lambda
}_{\operatorname{Im}}}\\
\vspace{-2.5mm} &  &  & \\
\frac{\partial G}{\partial\bar{v}_{\operatorname{Re}}} & \frac{\partial
G}{\partial\bar{v}_{\operatorname{Im}}} & \frac{\partial G}{\partial
\bar{\lambda}_{\operatorname{Re}}} & \frac{\partial G}{\partial\bar{\lambda
}_{\operatorname{Im}}}%
\end{bmatrix}
\in\mathbb{R}^{(2(n_{x}+1)n_{\xi})\times(2(n_{x}+1)n_{\xi})},
\]
where
\begin{align}
\frac{\partial F}{\partial\bar{v}_{\operatorname{Re}}}  &  =%
\begin{bmatrix}
\mathbb{E}[\Psi\Psi^{T}\!\otimes\!K]\!-\!\mathbb{E}[(\bar{\lambda
}_{\operatorname{Re}}{}^{T}\Psi)\Psi\Psi^{T}\!\otimes\!M_{\sigma}]\\
\vspace{-2.5mm}\\
-\mathbb{E}[(\bar{\lambda}_{\operatorname{Im}}{}^{T}\Psi)\Psi\Psi^{T}%
\!\otimes\!M_{\sigma}]
\end{bmatrix}
,\quad\frac{\partial F}{\partial\bar{\lambda}_{\operatorname{Re}}}=%
\begin{bmatrix}
-\mathbb{E}[\Psi^{T}\otimes(\Psi\Psi^{T}\!\otimes\!M_{\sigma})]\bar
{v}_{\operatorname{Re}}\\
\vspace{-2.5mm}\\
-\mathbb{E}[\Psi^{T}\otimes(\Psi\Psi^{T}\!\otimes\!M_{\sigma})]\bar
{v}_{\operatorname{Im}}%
\end{bmatrix}
,\label{eq:jac_Fu_1}\\
\frac{\partial F}{\partial\bar{v}_{\operatorname{Im}}}  &  =%
\begin{bmatrix}
\mathbb{E}[(\bar{\lambda}_{\operatorname{Im}}{}^{T}\Psi)\Psi\Psi^{T}%
\!\otimes\!M_{\sigma}]\\
\vspace{-2.5mm}\\
\mathbb{E}[\Psi\Psi^{T}\!\otimes\!K]\!-\!\mathbb{E}[(\bar{\lambda
}_{\operatorname{Re}}{}^{T}\Psi)\Psi\Psi^{T}\!\otimes\!M_{\sigma}]
\end{bmatrix}
,\quad\frac{\partial F}{\partial\bar{\lambda}_{\operatorname{Im}}}=%
\begin{bmatrix}
\mathbb{E}[\Psi^{T}\otimes(\Psi\Psi^{T}\!\otimes\!M_{\sigma})]\bar
{v}_{\operatorname{Im}}\\
\vspace{-2.5mm}\\
-\mathbb{E}[\Psi^{T}\otimes(\Psi\Psi^{T}\!\otimes\!M_{\sigma})]\bar
{v}_{\operatorname{Re}}%
\end{bmatrix}
, \label{eq:jac_Fu_2}%
\end{align}
and
\begin{align}
\frac{\partial G}{\partial\bar{v}_{\operatorname{Re}}}  &  =%
\begin{bmatrix}
2\mathbb{E}[\Psi\otimes(\bar{v}_{\operatorname{Re}}^{T}\Psi\Psi^{T}%
\!\otimes\!I_{n_{x}})]\\
\vspace{-2.5mm}\\
0
\end{bmatrix}
,\quad\quad\frac{\partial G}{\partial\bar{\lambda}_{\operatorname{Re}}%
}=0,\label{eq:jac_Gl_1}\\
\frac{\partial G}{\partial\bar{v}_{\operatorname{Im}}}  &  =%
\begin{bmatrix}
0\\
\vspace{-2.5mm}\\
2\mathbb{E}[\Psi\otimes(\bar{v}_{\operatorname{Im}}^{T}\Psi\Psi^{T}%
\!\otimes\!I_{n_{x}})]
\end{bmatrix}
,\quad\quad\frac{\partial G}{\partial\bar{\lambda}_{\operatorname{Im}}}=0.
\label{eq:jac_Gl_2}%
\end{align}
However, for convenience in the formulation of the preconditioners presented
later, we formulate Newton iteration with rescaled Jacobian
matrix$~\widehat{DJ}(\bar{v}_{\operatorname{Re}}^{(n)},\bar{v}%
_{\operatorname{Im}}^{(n)},\bar{\lambda}_{\operatorname{Re}}^{(n)}%
,\bar{\lambda}_{\operatorname{Im}}^{(n)})$, which at step$~n$ entails solving
a linear system
\begin{equation}%
\begin{bmatrix}
\frac{\partial F^{(n)}}{\partial\bar{v}_{\operatorname{Re}}} & \frac{\partial
F^{(n)}}{\partial\bar{v}_{\operatorname{Im}}} & \frac{\partial F^{(n)}%
}{\partial\bar{\lambda}_{\operatorname{Re}}} & \frac{\partial F^{(n)}%
}{\partial\bar{\lambda}_{\operatorname{Im}}}\\
\vspace{-2.5mm} &  &  & \\
-\frac{1}{2}\frac{\partial G^{(n)}}{\partial\bar{v}_{\operatorname{Re}}} &
-\frac{1}{2}\frac{\partial G^{(n)}}{\partial\bar{v}_{\operatorname{Im}}} & 0 &
0
\end{bmatrix}%
\begin{bmatrix}
\delta\bar{v}_{\operatorname{Re}}\\
\delta\bar{v}_{\operatorname{Im}}\\
\delta\bar{\lambda}_{\operatorname{Re}}\\
\delta\bar{\lambda}_{\operatorname{Im}}%
\end{bmatrix}
=-%
\begin{bmatrix}
F^{(n)}\\
-\frac{1}{2}G^{(n)}%
\end{bmatrix}
, \label{eq:Newton}%
\end{equation}
followed by an update
\begin{equation}%
\begin{bmatrix}
\bar{v}_{\operatorname{Re}}^{(n+1)}\\
\bar{v}_{\operatorname{Im}}^{(n+1)}\\
\bar{\lambda}_{\operatorname{Re}}^{(n+1)}\\
\bar{\lambda}_{\operatorname{Im}}^{(n+1)}%
\end{bmatrix}
=%
\begin{bmatrix}
\bar{v}_{\operatorname{Re}}^{(n)}\\
\bar{v}_{\operatorname{Im}}^{(n)}\\
\bar{\lambda}_{\operatorname{Re}}^{(n)}\\
\bar{\lambda}_{\operatorname{Im}}^{(n)}%
\end{bmatrix}
+%
\begin{bmatrix}
\delta\bar{v}_{\operatorname{Re}}\\
\delta\bar{v}_{\operatorname{Im}}\\
\delta\bar{\lambda}_{\operatorname{Re}}\\
\delta\bar{\lambda}_{\operatorname{Im}}%
\end{bmatrix}
. \label{eq:Newton-update}%
\end{equation}

\begin{remark}
We used the rescaling of the Jacobian in~\cite{Lee-2018-IMS} in order to
symmetrize the matrix in~(\ref{eq:Newton}), however we note that here it is
still in general nonsymmetric.
\end{remark}

Linear systems~(\ref{eq:Newton}) are solved inexactly using a preconditioned
GMRES method. We refer to Appendix~\ref{sec-app:ini} for the details of
evaluation of the right-hand side and of the matrix-vector product, and
to~\cite{Lee-2019-LRS} for a discussion of GMRES in a related context.

\subsection{Inexact Newton iteration}

As in~\cite{Lee-2018-IMS}, we consider a line-search modification of this
method following~\cite[Algorithm~11.4]{Nocedal-1999-NO} in order to improve
global convergence behavior of Newton iteration. Denoting
\[
\bar{v}^{(n)}=[(\bar{v}_{\operatorname{Re}}^{(n)})^{T},(\bar{v}%
_{\operatorname{Im}}^{(n)})^{T}]^{T}\quad\text{and}\quad\bar{\lambda}%
^{(n)}=[(\bar{\lambda}_{\operatorname{Re}}^{(n)})^{T},(\bar{\lambda
}_{\operatorname{Im}}^{(n)})^{T}]^{T},
\]
we define the merit function as the sum of squares,
\[
f(\bar{v}^{(n)},\bar{\lambda}^{(n)})=\frac{1}{2}\Vert{\widehat{r}}(\bar
{v}^{(n)},\bar{\lambda}^{(n)})\Vert_{2}^{2},\quad\text{where }\widehat{r}%
(\bar{v}^{(n)},\bar{\lambda}^{(n)})=%
\begin{bmatrix}
F^{(n)}\\
-\frac{1}{2}G^{(n)}%
\end{bmatrix}
,
\]
that is $\widehat{r}(\bar{v}^{(n)},\bar{\lambda}^{(n)})$ is the negative
right-hand side of~(\ref{eq:Newton}), i.e., it is a rescaled residual
of~(\ref{eq:nonlinear_res}), and we also denote
\[
f_{n}=f(\bar{v}^{(n)},\bar{\lambda}^{(n)}),\qquad\widehat{{r}}{_{n}%
=\widehat{r}(\bar{v}^{(n)},\bar{\lambda}^{(n)}),}\qquad\widehat{DJ}%
_{n}=\widehat{DJ}(\bar{v}^{(n)},\bar{\lambda}^{(n)}).
\]
As the initial approximation of the solution, we use the eigenvectors and
eigenvalues of the associated mean problem
\begin{equation}
K_{1}\left(  [v_{\operatorname{Re}}^{(0)}]_{1}+i[v_{\operatorname{Im}}%
^{(0)}]_{1}\right)  =\left(  [\lambda_{\operatorname{Re}}^{(0)}]_{1}%
+i[\lambda_{\operatorname{Im}}^{(0)}]_{1}\right)  M_{\sigma}\left(
[v_{\operatorname{Re}}^{(0)}]_{1}+i[v_{\operatorname{Im}}^{(0)}]_{1}\right)  ,
\label{eq:eigenproblem-mean}%
\end{equation}
concatenated by zeros, that is
\begin{align*}
\bar{v}^{(0)}  &  =\left[  ([v_{\operatorname{Re}}^{(0)}]_{1})^{T}%
,0,\ldots,([v_{\operatorname{Im}}^{(0)}]_{1})^{T},0,\ldots\right]  ^{T},\\
\bar{\lambda}^{(0)}  &  =\left[  [\lambda_{\operatorname{Re}}^{(0)}%
]_{1},0,\ldots,[\lambda_{\operatorname{Im}}^{(0)}]_{1},0,\ldots\right]  ^{T},
\end{align*}
and the initial residual is
\[
\widehat{r}_{0}=%
\begin{bmatrix}
F(\bar{v}^{(0)},\bar{\lambda}^{(0)})\\
-\frac{1}{2}G(\bar{v}^{(0)})
\end{bmatrix}
.
\]
The line-search Newton method is summarized in our setting as
Algorithm~\ref{alg:line_search}, and the choice of parameters$~\rho$ and$~c$
in the numerical experiments is discussed in Section~\ref{sec:numerical}.

\begin{algorithm}[hptb]
\caption{\cite[Algorithm 11.4]{Nocedal-1999-NO} Line-search Newton iteration}
\begin{algorithmic}[1]
\State Given $\rho,c \in (0, 1)$, set $\alpha^\ast = 1$.
\State Set $\bar v^{(0)}$ and $\bar \lambda^{(0)}$.
\For{$n=0,1,2,\ldots$}
\State $\widehat{DJ}_n p_n = - \widehat{r}_n$ \hfill (Solve inexactly to find the Newton update $p_n$.) \label{ln:lin_sys}
\State $\begin{bmatrix} \delta \bar v^{(n)} \\ \delta \bar \lambda^{(n)} \end{bmatrix} = p_n$
\State $\alpha_n = \alpha^\ast$
\While{ $f(\bar v^{(n)} + \alpha_n \delta \bar v^{(n)}, \bar \lambda^{(n)} + \alpha_n \delta \bar \lambda^{(n)}) > f_n + c \, \alpha_n \nabla f_n^T p_n  $}
\State $\alpha_n \leftarrow \rho \, \alpha_n$
\EndWhile
\State $\bar v^{(n+1)} \leftarrow \bar v^{(n)} + \alpha_n \delta \bar v^{(n)}$
\State $\bar \lambda^{(n+1)} \leftarrow \bar \lambda^{(n)} + \alpha_n \delta \bar \lambda^{(n)}$
\State Check for convergence. \label{ln:ls_conv-check}
\EndFor
\end{algorithmic}
\label{alg:line_search}
\end{algorithm}

The inexact iteration entails in each step a solution of the stochastic
Galerkin linear system in Line~\ref{ln:lin_sys} of
Algorithm~\ref{alg:line_search} given by~(\ref{eq:Newton}) using a Krylov
subspace method. We use preconditioned GMRES with the adaptive stopping
criteria,
\begin{equation}
\frac{\Vert\widehat{r}_{n}+\widehat{DJ}_{n}p_{n}\Vert_{2}}{\Vert
\widehat{r}_{n}\Vert_{2}}<\tau\left\Vert \widehat{r}{_{n-1}}\right\Vert _{2},
\label{eq:gmres-stop}%
\end{equation}
where $\tau=10^{-1}$. The for-loop is terminated when the convergence check in
Line~\ref{ln:ls_conv-check} is satisfied; in our numerical experiments we
check if $\left\Vert \widehat{r}_{n}\right\Vert _{2}<10^{-10}$.

Our implementation of the solvers\ is based on the so-called \emph{matricized}
formulation, in which we make use of isomorphism between~$%
\mathbb{R}
^{n_{x}n_{\xi}}$ and $%
\mathbb{R}
^{n_{x}\times n_{\xi}}$\ determined by the operators $\operatorname{vec}$ and
$\operatorname{mat}$: $\bar{v}=\operatorname{vec}(\bar{V})$, $\bar
{V}=\operatorname{mat}$($\bar{v})$, where $\bar{v}\in%
\mathbb{R}
^{n_{x}n_{\xi}}$, $\bar{V}\in%
\mathbb{R}
^{n_{x}\times n_{\xi}}$. The upper/lower case notation is assumed throughout
the paper, so $\bar{R}=\operatorname{mat}$($\bar{r})$, etc. Specifically, we
define the \emph{matricized}\ coefficients of the eigenvector expansion
\begin{equation}
\bar{V}=\operatorname{mat}(\bar{v})=\left[  v_{1},v_{2},\ldots,v_{n_{\xi}%
}\right]  \in%
\mathbb{R}
^{n_{x}\times n_{\xi}}, \label{eq:U}%
\end{equation}
where the column~$k$ contains the coefficients associated with the basis
function$~\psi_{k}$. A detailed formulation of the GMRES in the matricized
format can be found, e.g., in~\cite{Lee-2019-LRS}. We only describe
computation of the matrix-vector product (Appendix~\ref{sec-app:ini}), and in
the next section we formulate several preconditioners.

\subsection{Preconditioners for the Newton iteration}

\label{sec:preconditioners} In order to allow for an efficient iterative
solution of linear systems in Line~\ref{ln:lin_sys} of
Algorithm~\ref{alg:line_search} given by~(\ref{eq:Newton}) using a Krylov
subspace method, it is necessary to provide a preconditioner. In this section,
we adapt the mean-based preconditioner and two of the constraint
preconditioners from~\cite{Lee-2018-IMS} to the formulation with separated
real and complex parts, and we write them in the matricized format. The idea
can be motivated as follows. The preconditioners are based on approximations
of the blocks in~(\ref{eq:jac-scheme}). The mean-based preconditioner is
\textcolor{black}{inspired by} the approximation
\[
\left[
\begin{array}
[c]{cc}%
\overline{A} & 0\\
0 & \overline{S}%
\end{array}
\right]  ,
\]
where $\overline{A}\approx\left[  A_{\operatorname{Re}}\;A_{\operatorname{Im}%
}\right]  $ and the Schur complement $\overline{S}\approx-\frac{1}{2}\left[
C_{\operatorname{Re}}\;C_{\operatorname{Im}}\right]  \left[
A_{\operatorname{Re}}\;A_{\operatorname{Im}}\right]  ^{-1}\left[
B_{\operatorname{Re}}\;B_{\operatorname{Im}}\right]  $. The constraint
preconditioners are based on the approximation
\[
\left[
\begin{array}
[c]{cc}%
\overline{A} & \overline{B}\\
\overline{C} & 0
\end{array}
\right]  ,
\]
where $\overline{B}\approx\left[  B_{\operatorname{Re}}\;B_{\operatorname{Im}%
}\right]  $ and $\overline{C}\approx-\frac{1}{2}\left[  C_{\operatorname{Re}%
}\;C_{\operatorname{Im}}\right]  $. Next, considering the truncation of the
series in~(\ref{eq:jac_A_Re})--(\ref{eq:jac_C_Im}) to the very first term, we
get
\begin{align*}
\overline{A}  &  \approx I_{n_{\xi}}\otimes\left[
\begin{array}
[c]{cc}%
K_{1}-\lambda_{\operatorname{Re},1}M_{\sigma} & \lambda_{\operatorname{Im}%
,1}M_{\sigma}\\
-\lambda_{\operatorname{Im},1}M_{\sigma} & K_{1}-\lambda_{\operatorname{Re}%
,1}M_{\sigma}%
\end{array}
\right]  \text{ see left columns in (\ref{eq:jac_Fu_1})--(\ref{eq:jac_Fu_2})
and (\ref{eq:jac_A_Re})--(\ref{eq:jac_A_Im})},\\
\overline{B}  &  \approx I_{n_{\xi}}\otimes\left[
\begin{array}
[c]{cc}%
-M_{\sigma}v_{\operatorname{Re},1} & M_{\sigma}v_{\operatorname{Im},1}\\
-M_{\sigma}v_{\operatorname{Im},1} & -M_{\sigma}v_{\operatorname{Re},1}%
\end{array}
\right]  \text{ see right columns in~(\ref{eq:jac_Fu_1})--(\ref{eq:jac_Fu_2})
and (\ref{eq:jac_B_Re})--(\ref{eq:jac_B_Im}),}\\
\overline{C}  &  \approx I_{n_{\xi}}\otimes\left[
\begin{array}
[c]{cc}%
-v_{\operatorname{Re},1}^{T} & 0_{n_{x}\times1}\\
0_{n_{x}\times1} & -v_{\operatorname{Im},1}^{T}%
\end{array}
\right]  \text{ see (\ref{eq:jac_Gl_1})--(\ref{eq:jac_Gl_2}) and
(\ref{eq:jac_C_Re})--(\ref{eq:jac_C_Im}).}%
\end{align*}
The precise formulations are listed for the mean-based preconditioner as
Algorithm~\ref{alg:MB} and for the constraint mean-based preconditioner as
Algorithm~\ref{alg:cMB}. Finally, the constraint hierarchical Gauss-Seidel
preconditioner is listed as Algorithm~\ref{alg:chGS}. It can be viewed as an
extension of Algorithm~\ref{alg:cMB}, because the solves with stochastic
Galerkin matrices~(\ref{eq:M-preconditioner}) are used also in this
preconditioner, but in addition the right-hand sides are updated using an idea
inspired by Gauss-Seidel method in a for-loop over the degree of the
gPC\ basis. Moreover, as proposed in~\cite{Sousedik-2014-THP}, the
matrix-vector multiplications, used in the setup of the right-hand sides can
be truncated in the sense that in the summations, $t=1,\dots,n_{\xi}$ is
replaced by $t\in\mathcal{I}_{t}$, where $\mathcal{I}_{t}=\left\{
1,\dots,n_{t}\right\}  $ with $n_{t}=\binom{m_{\xi}+p_{t}}{p_{t}}$ for some
$p_{t}\leq p$, and in particular with $p_{t}=0$ the chGS preconditioner
(Algorithm~\ref{alg:chGS}) reduces to the cMB{\ preconditioner
(Algorithm~\ref{alg:cMB}). We also note that, since the initial guess is zero
in }Algorithm~\ref{alg:chGS}, the multiplications by$~\mathcal{F}_{1}$ and
$\mathcal{F}_{d+1}$\ vanish from (\ref{eq:algchGS1})--(\ref{eq:algchGS2}).

\begin{algorithm}[hptb]
\caption{Mean-based preconditioner (MB)}
\label{alg:MB}
The preconditioner $\mathcal{M}_{\text{MB}}:\left(  \bar{R}%
^{v_{\operatorname{Re}}},\bar{R}^{v_{\operatorname{Im}}},\bar{R}%
^{\lambda_{\operatorname{Re}}},\bar{R}^{\lambda_{\operatorname{Im}}}\right)
\longmapsto\left(  \bar{V}^{v_{\operatorname{Re}}},\bar{V}%
^{v_{\operatorname{Im}}},\bar{V}^{\lambda_{\operatorname{Re}}},\bar
{V}^{\lambda_{\operatorname{Im}}}\right)  $\ is defined as
\begin{equation}
\mathcal{M}_{\text{MB}}\left[
\begin{array}
[c]{c}%
\bar{V}^{v_{\operatorname{Re}}}\\
\bar{V}^{v_{\operatorname{Im}}}\\
\bar{V}^{\lambda_{\operatorname{Re}}}\\
\bar{V}^{\lambda_{\operatorname{Im}}}%
\end{array}
\right]  =\left[
\begin{array}
[c]{c}%
\bar{R}^{v_{\operatorname{Re}}}\\
\bar{R}^{v_{\operatorname{Im}}}\\
\bar{R}^{\lambda_{\operatorname{Re}}}\\
\bar{R}^{\lambda_{\operatorname{Im}}}%
\end{array}
\right]  .\label{eq:algMB}%
\end{equation}
Above
\begin{equation}
\mathcal{M}_{\text{MB}}=\left[
\begin{array}
[c]{cc}%
\widetilde{A} & 0_{2n_{x}\times2}\\
0_{2\times2n_{x}} & \left[
\begin{array}
[c]{cc}%
-w_{\operatorname{Re}}^{(0)T} & 0_{n_{x}\times1}\\
0_{n_{x}\times1} & -w_{\operatorname{Im}}^{(0)T}%
\end{array}
\right]  \widetilde{A}^{-1} \left[
\begin{array}
[c]{cc}%
- M_{\sigma} w_{\operatorname{Re}}^{(0)} & M_{\sigma} w_{\operatorname{Im}}^{(0)}\\
- M_{\sigma} w_{\operatorname{Im}}^{(0)} & - M_{\sigma} w_{\operatorname{Re}}^{(0)}%
\end{array}
\right]
\end{array}
\right]  ,\label{eq:MB-preconditioner}%
\end{equation}
where~$w_{\operatorname{Re}}^{(0)}$ and $w_{\operatorname{Im}}^{(0)}$ are 
the real and imaginary\ components of eigenvector $w$ of the stencil
$(K_{1},M_{\sigma})$ with corresponding eigenvalue $\mu=\mu_{\operatorname{Re}%
}$+i$\mu_{\operatorname{Im}}$, cf.~(\ref{eq:eigenproblem-mean}), 
and
\begin{equation}
\widetilde{A}=\left[
\begin{array}
[c]{cc}%
K_{1}-\epsilon_{\operatorname{Re}}\mu_{\operatorname{Re}}M_{\sigma} &
\epsilon_{\operatorname{Im}}\mu_{\operatorname{Im}}M_{\sigma}\\
-\epsilon_{\operatorname{Im}}\mu_{\operatorname{Im}}M_{\sigma} &
K_{1}-\epsilon_{\operatorname{Re}}\mu_{\operatorname{Re}}M_{\sigma}%
\end{array}
\right]  ,\label{eq:tildeA}%
\end{equation}
with constants $\epsilon_{\operatorname{Re}}$, $\epsilon_{\operatorname{Im}}$
further specified in the numerical experiments section.
\end{algorithm}

\begin{algorithm}[hptb]
\caption{Constraint mean-based preconditioner (cMB)}
\label{alg:cMB}
The preconditioner $\mathcal{M}_{\text{cMB}}:\left(  \bar{R}%
^{v_{\operatorname{Re}}},\bar{R}^{v_{\operatorname{Im}}},\bar{R}%
^{\lambda_{\operatorname{Re}}},\bar{R}^{\lambda_{\operatorname{Im}}}\right)
\longmapsto\left(  \bar{V}^{v_{\operatorname{Re}}},\bar{V}%
^{v_{\operatorname{Im}}},\bar{V}^{\lambda_{\operatorname{Re}}},\bar
{V}^{\lambda_{\operatorname{Im}}}\right)  $\ is defined as
\begin{equation}
\mathcal{M}_{\text{cMB}}\left[
\begin{array}
[c]{c}%
\bar{V}^{v_{\operatorname{Re}}}\\
\bar{V}^{v_{\operatorname{Im}}}\\
\bar{V}^{\lambda_{\operatorname{Re}}}\\
\bar{V}^{\lambda_{\operatorname{Im}}}%
\end{array}
\right]  =\left[
\begin{array}
[c]{c}%
\bar{R}^{v_{\operatorname{Re}}}\\
\bar{R}^{v_{\operatorname{Im}}}\\
\bar{R}^{\lambda_{\operatorname{Re}}}\\
\bar{R}^{\lambda_{\operatorname{Im}}}%
\end{array}
\right]  .\label{eq:algcMB}%
\end{equation}
Above
\begin{equation}
\mathcal{M}_{\text{cMB}}=\left[
\begin{array}
[c]{cc}%
\widetilde{A} &
\begin{array}
[c]{cc}%
- M_{\sigma} w_{\operatorname{Re}}^{(0)} & M_{\sigma} w_{\operatorname{Im}}^{(0)}\\
- M_{\sigma} w_{\operatorname{Im}}^{(0)} & - M_{\sigma} w_{\operatorname{Re}}^{(0)}%
\end{array}
\\%
\begin{array}
[c]{cc}%
-w_{\operatorname{Re}}^{(0)T} & 0_{n_{x}\times1}\\
0_{n_{x}\times1} & -w_{\operatorname{Im}}^{(0)T}%
\end{array}
& 0_{2\times2}%
\end{array}
\right]  ,\label{eq:M-preconditioner}%
\end{equation}
with $w_{\operatorname{Re}}^{(0)}$, $w_{\operatorname{Im}}^{(0)}$ and $\widetilde{A}$ set as in Algorithm~\ref{alg:MB}.
\end{algorithm}

\begin{algorithm}[hptb]
\caption{Constraint hierarchical Gauss-Seidel preconditioner (chGS)}
\label{alg:chGS}
The preconditioner $M_{\text{chGS}}:\left(  \bar{R}%
^{v_{\operatorname{Re}}},\bar{R}^{v_{\operatorname{Im}}},\bar{R}%
^{\lambda_{\operatorname{Re}}},\bar{R}^{\lambda_{\operatorname{Im}}}\right)
\longmapsto\left(  \bar{V}^{v_{\operatorname{Re}}},\bar{V}%
^{v_{\operatorname{Im}}},\bar{V}^{\lambda_{\operatorname{Re}}},\bar
{V}^{\lambda_{\operatorname{Im}}}\right)  $\
is defined as follows.
\begin{algorithmic}[1]
\State Set the initial solution $\left( \bar{V}^{v\RE}, \bar{V}^{v\IM},\bar{V}^{\lambda\RE
},\bar{V}^{\lambda\IM}\right) $ to zero and update as: 
\State Solve
\vspace{-5mm}
\begin{equation}
\mathcal{M}_{\text{cMB}}
\left(
\begin{array}{c}
V_{1}^{v\RE} \\
V_{1}^{v\IM}\\
V_{1}^{\lambda\RE}\\
V_{1}^{\lambda\IM}
\end{array}%
\right) =\left[
\begin{array}{c}
R_{1}^{v\RE} \\
R_{1}^{v\IM} \\
R_{1}^{\lambda\RE}\\
R_{1}^{\lambda\IM}%
\end{array}%
\right] -\mathcal{F}_{1}\left(
\begin{array}{c}
V_{\left( 2:n_{\xi }\right) }^{v\RE} \\
V_{\left( 2:n_{\xi }\right) }^{v\IM} \\
V_{\left( 2:n_{\xi }\right) }^{\lambda\RE}\\
V_{\left( 2:n_{\xi }\right) }^{\lambda\IM}
\end{array}%
\right),   \label{eq:algchGS1}
\end{equation}%
where $\mathcal{M}_{\text{cMB}}$ is set as in Algorithm~\ref{alg:cMB}, and
\[
\mathcal{F}_{1}\left(
\begin{array}
[c]{c}%
V_{1}^{v\RE}\\
\vspace{-3mm}\\
V_{1}^{v\IM}\\
\vspace{-3mm}\\
V_{1}^{\lambda\RE}\\
\vspace{-3mm}\\
V_{1}^{\lambda\IM}%
\end{array}
\right)  =\!\sum_{t\in\mathcal{I}_{t}}\!\left[  \!\!%
\begin{array}
[c]{c}%
\mathcal{F}_{1}^{A}\left[  h_{t,\left(  2:n_{\xi}\right)  \left(  1\right)
}\right]  \\
\vspace{-3mm}\\
\mathcal{F}_{1}^{B}\left[  h_{t,\left(  2:n_{\xi}\right)  \left(  1\right)
}\right]  \\
\vspace{-3mm}\\
-(v_{\text{Re},t}^{(n)})^{T}V_{\left(  2:n_{\xi}\right)  }^{v\RE}\left[
h_{t,\left(  2:n_{\xi}\right)  \left(  1\right)  }\right]  \\
\vspace{-3mm}\\
-(v_{\text{Im},t}^{(n)})^{T}V_{\left(  2:n_{\xi}\right)  }^{v\IM}\left[
h_{t,\left(  2:n_{\xi}\right)  \left(  1\right)  }\right]
\end{array}
\!\!\right]  ,
\]
\vspace{-3mm}
\begin{align*}
\mathcal{F}_{1}^{A}  & =\left(  \left(  \!K_{t}\!-\!\lambda_{\text{Re}%
,t}^{(n)}M_{\sigma}\!\right)  \!V_{\left(  2:n_{\xi}\right)  }^{v\RE}%
\!+\!\lambda_{\text{Im},t}^{(n)}M_{\sigma}\!V_{\left(  2:n_{\xi}\right)
}^{v\IM}\!\!-\!v_{\text{Re},t}^{(n)}M_{\sigma}\!V_{\left(  2:n_{\xi}\right)
}^{\lambda\RE}\!+\!v_{\text{Im},t}^{(n)}M_{\sigma}\!V_{\left(  2:n_{\xi
}\right)  }^{\lambda\IM}\!\right)  ,\\
\mathcal{F}_{1}^{B}  & =\left(  \left(  \!K_{t}\!-\!\lambda_{\text{Re}%
,t}^{(n)}M_{\sigma}\!\right)  \!V_{\left(  2:n_{\xi}\right)  }^{v\IM}%
\!-\!\lambda_{\text{Im},t}^{(n)}M_{\sigma}\!V_{\left(  2:n_{\xi}\right)
}^{v\RE}\!\!-\!v_{\text{Re},t}^{(n)}M_{\sigma}\!V_{\left(  2:n_{\xi}\right)
}^{\lambda\IM}\!-\!v_{\text{Im},t}^{(n)}M_{\sigma}\!V_{\left(  2:n_{\xi
}\right)  }^{\lambda\RE}\!\right)  ,
\end{align*}
and $v_{\text{Re},t}^{(n)}$, $v_{\text{Im},t}^{(n)}$ the $t$th gPC
coefficients of eigenvector $v^{(n)}$ at step $n$ of
Algorithm~\ref{alg:line_search}.
\For{$d=1,\ldots, p-1$}
\State Set $\ell =\left( n_{\ell }+1:n_{u}\right) ,\text{ where }n_{\ell }=\binom{n_{\xi
}+d-1}{d-1}\text{ and }n_{u}=\binom{n_{\xi }+d}{d}$.
\State Solve
\vspace{-5mm}
\begin{equation}
\mathcal{M}_{\text{cMB}}
\left(
\begin{array}{c}
V_{\left( \ell \right) }^{v\RE} \\
V_{\left( \ell \right) }^{v\IM} \\
V_{\left( \ell \right) }^{\lambda\RE}\\
V_{\left( \ell \right) }^{\lambda\IM}%
\end{array}%
\right) =\left[
\begin{array}{c}
R_{\left( \ell \right)}^{v\RE} \\
R_{\left( \ell \right)}^{v\IM} \\
R_{\left( \ell \right) }^{\lambda\RE}\\
R_{\left( \ell \right) }^{\lambda\IM}%
\end{array}%
\right] -\mathcal{E}_{d+1}\left(
\begin{array}{c}
V_{(1:n_{\ell })}^{v\RE} \\
V_{(1:n_{\ell })}^{v\IM} \\
V_{(1:n_{\ell })}^{\lambda\RE}\\
V_{(1:n_{\ell })}^{\lambda\IM}%
\end{array}%
\right) -\mathcal{F}_{d+1}\left(
\begin{array}{c}
V_{(n_{u}+1:n_{\xi })}^{v\RE} \\
V_{(n_{u}+1:n_{\xi })}^{v\IM} \\
V_{(n_{u}+1:n_{\xi })}^{\lambda\RE}\\
V_{(n_{u}+1:n_{\xi })}^{\lambda\IM}
\end{array}%
\right) ,
\label{eq:algchGS2}
\end{equation}
where
\vspace{-5mm}
\begin{align*}
\mathcal{E}_{d+1}\left(  \!\!\!%
\begin{array}
[c]{c}%
V_{(1:n_{\ell})}^{v\RE}\\
V_{(1:n_{\ell})}^{v\IM}\\
V_{(1:n_{\ell})}^{\lambda\RE}\\
V_{(1:n_{\ell})}^{\lambda\IM}%
\end{array}
\!\!\!\right)  \! &  =\!\sum_{t\in\mathcal{I}_{t}}\!\!\left[  \!\!%
\begin{array}
[c]{c}%
\mathcal{E}_{d+1}^{A}\left[  h_{t,\left(  1:n_{\ell}\right)  \left(
\ell\right)  }\right]  \\
\vspace{-3mm}\\
\mathcal{E}_{d+1}^{B}\left[  h_{t,\left(  1:n_{\ell}\right)  \left(
\ell\right)  }\right]  \\
\vspace{-3mm}\\
-(v_{\text{Re},t}^{(n)})^{T}V_{\left(  1:n_{\ell}\right)  }^{v\RE}\left[
h_{t,\left(  1:n_{\ell}\right)  \left(  \ell\right)  }\right]  \\
\vspace{-3mm}\\
-(v_{\text{Im},t}^{(n)})^{T}V_{\left(  1:n_{\ell}\right)  }^{v\IM}\left[
h_{t,\left(  1:n_{\ell}\right)  \left(  \ell\right)  }\right]
\end{array}
\!\!\right]  ,\\
\mathcal{F}_{d+1}\left(  \!\!\!%
\begin{array}
[c]{c}%
V_{(n_{u}+1:n_{\xi})}^{v\RE}\\
V_{(n_{u}+1:n_{\xi})}^{v\IM}\\
V_{(n_{u}+1:n_{\xi})}^{\lambda\RE}\\
V_{(n_{u}+1:n_{\xi})}^{\lambda\IM}%
\end{array}
\!\!\!\right)  \! &  =\!\sum_{t\in\mathcal{I}_{t}}\!\!\left[  \!%
\begin{array}
[c]{c}%
\mathcal{F}_{d+1}^{A}\left[  h_{t,\left(  n_{u}+1:n_{\xi}\right)  \left(
\ell\right)  }\right]  \\
\vspace{-3mm}\\
\mathcal{F}_{d+1}^{B}\left[  h_{t,\left(  n_{u}+1:n_{\xi}\right)  \left(
\ell\right)  }\right]  \\
\vspace{-3mm}\\
-(v_{\text{Re},t}^{(n)})^{T}V_{\left(  n_{u}+1:n_{\xi}\right)  }^{v\RE}\left[
h_{t,\left(  n_{u}+1:n_{\xi}\right)  \left(  \ell\right)  }\right]  \\
\vspace{-3mm}\\
-(v_{\text{Im},t}^{(n)})^{T}V_{\left(  n_{u}+1:n_{\xi}\right)  }^{v\IM}\left[
h_{t,\left(  n_{u}+1:n_{\xi}\right)  \left(  \ell\right)  }\right]
\end{array}
\!\!\right]  ,
\end{align*}
\vspace{-3mm}
\begin{align*}
\mathcal{E}_{d+1}^{A} &  =\left(  \left(  \!K_{t}\!-\!\lambda_{\text{Re}%
,t}^{(n)}M_{\sigma}\!\right)  \!V_{\left(  1:n_{\ell}\right)  }^{v\RE}%
\!+\!\lambda_{\text{Im},t}^{(n)}M_{\sigma}\!V_{\left(  1:n_{\ell}\right)
}^{v\IM}\!\!-\!v_{\text{Re},t}^{(n)}M_{\sigma}\!V_{\left(  1:n_{\ell}\right)
}^{\lambda\RE}\!+\!v_{\text{Im},t}^{(n)}M_{\sigma}\!V_{\left(  1:n_{\ell
}\right)  }^{\lambda\IM}\!\right)  ,\\
\mathcal{E}_{d+1}^{B} &  =\left(  \left(  \!K_{t}\!-\!\lambda_{\text{Re}%
,t}^{(n)}M_{\sigma}\!\right)  \!V_{\left(  1:n_{\ell}\right)  }^{v\IM}%
\!-\!\lambda_{\text{Im},t}^{(n)}M_{\sigma}\!V_{\left(  1:n_{\ell}\right)
}^{v\RE}\!\!-\!v_{\text{Re},t}^{(n)}M_{\sigma}\!V_{\left(  1:n_{\ell}\right)
}^{\lambda\IM}\!-\!v_{\text{Im},t}^{(n)}M_{\sigma}\!V_{\left(  1:n_{\ell
}\right)  }^{\lambda\RE}\!\right)  ,\\
\mathcal{F}_{d+1}^{A} &  =\left(  \left(  \!K_{t}\!-\!\lambda_{\text{Re}%
,t}^{(n)}M_{\sigma}\!\right)  \!V_{\left(  n_{u}+1:n_{\xi}\right)  }%
^{v\RE}\!+\!\lambda_{\text{Im},t}^{(n)}M_{\sigma}\!V_{\left(  n_{u}+1:n_{\xi
}\right)  }^{v\IM}\!\!-\!v_{\text{Re},t}^{(n)}M_{\sigma}\!V_{\left(
n_{u}+1:n_{\xi}\right)  }^{\lambda\RE}\!+\!v_{\text{Im},t}^{(n)}M_{\sigma
}\!V_{\left(  n_{u}+1:n_{\xi}\right)  }^{\lambda\IM}\!\right)  ,\\
\mathcal{F}_{d+1}^{B} &  =\left(  \left(  \!K_{t}\!-\!\lambda_{\text{Re}%
,t}^{(n)}M_{\sigma}\!\right)  \!V_{\left(  n_{u}+1:n_{\xi}\right)  }%
^{v\IM}\!-\!\lambda_{\text{Im},t}^{(n)}M_{\sigma}\!V_{\left(  n_{u}+1:n_{\xi
}\right)  }^{v\RE}\!\!-\!v_{\text{Re},t}^{(n)}M_{\sigma}\!V_{\left(
n_{u}+1:n_{\xi}\right)  }^{\lambda\IM}\!-\!v_{\text{Im},t}^{(n)}M_{\sigma
}\!V_{\left(  n_{u}+1:n_{\xi}\right)  }^{\lambda\RE}\!\right)  .
\end{align*}
\EndFor
\algstore{split_here}
\end{algorithmic}
\end{algorithm}

\begin{algorithm}[hptb]
\caption{Constraint hierarchical Gauss-Seidel preconditioner (chGS), continued}
\label{alg:chGS_cont}
\begin{algorithmic}[1]
\algrestore{split_here}
\State Set $\ell =\left( n_{u}+1:n_{\xi }\right)$.
\State
Solve%
\vspace{-5mm}
\[
\mathcal{M}_{\text{cMB}}\left(
\begin{array}
[c]{c}%
V_{(\ell)}^{v\RE}\\
V_{(\ell)}^{v\IM}\\
V_{(\ell)}^{\lambda\RE}\\
V_{(\ell)}^{\lambda\IM}%
\end{array}
\right)  =\left[
\begin{array}
[c]{c}%
R_{(\ell)}^{v\RE}\\
R_{(\ell)}^{v\IM}\\
R_{(\ell)}^{\lambda\RE}\\
R_{(\ell)}^{\lambda\IM}%
\end{array}
\right]  -\mathcal{E}_{p+1}\left(
\begin{array}
[c]{c}%
V_{(1:n_{u})}^{v\RE}\\
V_{(1:n_{u})}^{v\IM}\\
V_{(1:n_{u})}^{\lambda\RE}\\
V_{(1:n_{u})}^{\lambda\IM}%
\end{array}
\right)  ,
\]
where
\vspace{-3mm}
\[
\mathcal{E}_{p+1}\left(
\begin{array}
[c]{c}%
V_{(1:n_{u})}^{v\RE}\\
V_{(1:n_{u})}^{v\IM}\\
V_{(1:n_{u})}^{\lambda\RE}\\
V_{(1:n_{u})}^{\lambda\IM}%
\end{array}
\right)  =\sum_{t\in\mathcal{I}_{t}}\left[  \!\!%
\begin{array}
[c]{c}%
\mathcal{E}_{p+1}^{A}\left[  h_{t,\left(  1:n_{u}\right)  \left(  \ell\right)
}\right] \\
\vspace{-3mm}\\
\mathcal{E}_{p+1}^{B}\left[  h_{t,\left(  1:n_{u}\right)  \left(  \ell\right)
}\right] \\
\vspace{-3mm}\\
-(v_{\text{Re},t}^{(n)})^{T}V_{\left(  1:n_{u}\right)  }^{v\RE}\left[
h_{t,\left(  1:n_{u}\right)  \left(  \ell\right)  }\right] \\
\vspace{-3mm}\\
-(v_{\text{Im},t}^{(n)})^{T}V_{\left(  1:n_{u}\right)  }^{v\IM}\left[
h_{t,\left(  1:n_{u}\right)  \left(  \ell\right)  }\right]
\end{array}
\!\!\right]  ,
\]
\vspace{-3mm}
\begin{align*}
\mathcal{E}_{p+1}^{A}  &  =\left(  \left(  \!K_{t}\!-\!\lambda_{\text{Re}%
,t}^{(n)}M_{\sigma}\!\right)  \!V_{\left(  1:n_{u}\right)  }^{v\RE}%
\!+\!\lambda_{\text{Im},t}^{(n)}M_{\sigma}\!V_{\left(  1:n_{u}\right)
}^{v\IM}\!\!-\!v_{\text{Re},t}^{(n)}M_{\sigma}\!V_{\left(  1:n_{u}\right)
}^{\lambda\RE}\!+\!v_{\text{Im},t}^{(n)}M_{\sigma}\!V_{\left(  1:n_{u}\right)
}^{\lambda\IM}\!\right)  ,\\
\mathcal{E}_{p+1}^{B}  &  =\left(  \left(  \!K_{t}\!-\!\lambda_{\text{Re}%
,t}^{(n)}M_{\sigma}\!\right)  \!V_{\left(  1:n_{u}\right)  }^{v\IM}%
\!-\!\lambda_{\text{Im},t}^{(n)}M_{\sigma}\!V_{\left(  1:n_{u}\right)
}^{v\RE}\!\!-\!v_{\text{Re},t}^{(n)}M_{\sigma}\!V_{\left(  1:n_{u}\right)
}^{\lambda\IM}\!-\!v_{\text{Im},t}^{(n)}M_{\sigma}\!V_{\left(  1:n_{u}\right)
}^{\lambda\RE}\!\right)  .
\end{align*}
\For{$d=p-1,\ldots, 1$}
\State Set $\ell =\left( n_{\ell }+1:n_{u}\right) ,\text{ where }n_{\ell }=\binom{n_{\xi
}+d-1}{d-1}\text{ and }n_{u}=\binom{n_{\xi }+d}{d}$.
\State  Solve~(\ref{eq:algchGS2}).
\EndFor
\State Solve~(\ref{eq:algchGS1}).
\end{algorithmic}
\end{algorithm}

\subsubsection{Updating the constraint preconditioner}

\label{sec:SMW} It is also possible to update the application of the
constraint mean-based preconditioner in order to incorporate the latest
approximations of the eigenvector mean coefficients represented by the
vectors$~w_{\operatorname{Re}/\operatorname{Im}}^{(n)}$
in~(\ref{eq:M-preconditioner}). Suppose the inverse of the matrix$~\mathcal{M}%
_{\text{cMB}}$ from~(\ref{eq:M-preconditioner}) for $n=0$ is available, e.g.,
as the LU decomposition. That is, we have the preconditioner for the initial
step of the Newton iteration, and let us denote its application by $X^{-1}$.
Specifically, $X^{-1}=U^{-1}L^{-1}$, where $L$ and $U$ are the factors
of$~\mathcal{M}_{\text{cMB}}$. Next, let us consider two matrices $Y$ and $Z$
such that
\begin{equation}
YZ^{T}=\left[
\begin{array}
[c]{cc}%
0_{2n_{x}\times2n_{x}} &
\begin{array}
[c]{cc}%
-M_{\sigma}\Delta w_{\operatorname{Re}}^{(n)} & M_{\sigma}\Delta
w_{\operatorname{Im}}^{(n)}\\
-M_{\sigma}\Delta w_{\operatorname{Im}}^{(n)} & -M_{\sigma}\Delta
w_{\operatorname{Re}}^{(n)}%
\end{array}
\\%
\begin{array}
[c]{cc}%
-\Delta w_{\operatorname{Re}}^{(n)} & 0_{n_{x}\times1}\\
0_{n_{x}\times1} & -\Delta w_{\operatorname{Im}}^{(n)}%
\end{array}
& 0_{2\times2}%
\end{array}
\right]  , \label{eq:M-preconditione-update}%
\end{equation}
where
\[
\Delta w_{\operatorname{Re}}^{(n)}=w_{\operatorname{Re}}^{(n)}%
-w_{\operatorname{Re}}^{(0)},\qquad\Delta w_{\operatorname{Im}}^{(n)}%
=w_{\operatorname{Im}}^{(n)}-w_{\operatorname{Im}}^{(0)}.
\]
Specifically, $YZ^{T}$ is the rank$~4$ update of the preconditioner at
step$~n$ of Newton iteration, and the matrices $Y$ and $Z$ can be computed
using a sparse SVD,
\textcolor{black}{which is available, e.g., in~\cite{MATLAB:R2021a_u3}}. In
implementation, using \textsc{Matlab} notation with $YZ^{T}%
=\mathtt{\mathtt{YZt}}$, we compute $\mathtt{[U,S,V]=svds(\mathtt{YZt,4})}$
and set $Y=\mathtt{U\ast S}$, $Z^{T}=\mathtt{\mathtt{V}}^{\prime}$. Finally, a
solve $\mathcal{M}_{\text{cMB}}v=u$ at step$~n$ of Newton iteration may be
facilitated using the Sherman-Morrison-Woodbury formula see,
e.g.~\cite{Hager-1989-UIM}, or \cite[Section~3.8]{Meyer-2000-MAA}, as
\[
v=\left(  X+YZ^{T}\right)  ^{-1}u=\left(  X^{-1}-X^{-1}Y\left(  I+Z^{T}%
X^{-1}Y\right)  ^{-1}Z^{T}X^{-1}\right)  u.
\]

\section{Bifurcation analysis of Navier\textendash Stokes equations with
stochastic viscosity}

\label{sec:NS}
Here, we follow the setup from~\cite{Sousedik-2016-SGM} and assume that the
viscosity$~\nu$ is given by a stochastic\ expansion
\begin{equation}
\nu(x,\xi)=\sum_{\ell=1}^{n_{\nu}}\nu_{\ell}(x)\,\psi_{\ell}(\xi),
\label{eq:viscosity}%
\end{equation}
where $\left\{  \nu_{\ell}(x)\right\}  $\ is a set of given deterministic
spatial functions. We note that there are several possible interpretations of
such setup~\cite{Powell-2012-PSS,Sousedik-2016-SGM}. Assuming fixed geometry,
the stochastic viscosity is equivalent to Reynolds number being stochastic
and, for example, to a scenario when the volume of fluid moving into a channel
is uncertain. Consider the discretization of~(\ref{eq:NS-time}) by a
div-stable mixed finite element method,
and let the bases for the velocity and pressure spaces be denoted $\left\{
\phi_{i}\right\}  _{i=1}^{n_{u}}$ and $\left\{  \varphi_{j}\right\}
_{i=1}^{n_{p}}$, respectively. We further assume that we have a discrete
approximation of the steady-state solution of the stochastic counterpart
of~(\ref{eq:NS-time}), given as\footnote{Techniques for computing these
approximations were studied in~\cite{Lee-2019-LRS,Sousedik-2016-SGM}.}
\begin{align*}
\vec{u}(x,\xi)  &  {\color{black}\approx} \sum_{k=1}^{n_{\xi}}\sum_{i=1}%
^{n_{u}}u_{ik}\phi_{i}(x)\psi_{k}(\xi)=\sum_{k=1}^{n_{\xi}}\vec{u}_{k}(x)\psi
_{k}(\xi),\\
p(x,\xi)  &  {\color{black}\approx} \sum_{k=1}^{n_{\xi}}\sum_{j=1}^{n_{p}%
}p_{jk}\varphi_{j}(x)\psi_{k}(\xi)=\sum_{k=1}^{n_{\xi}}p_{k}(x)\psi_{k}(\xi).
\end{align*}

We are interested in a stochastic counterpart of the generalized eigenvalue
problem~(\ref{eq:eig-det}), which we write as
\begin{equation}
\mathcal{J(\xi)}v=\lambda M_{\sigma}v, \label{eq:NS-eig-gen}%
\end{equation}
where$~M_{\sigma}$ is the deterministic (shifted) mass matrix
from~(\ref{eq:M-shifted}), and $\mathcal{J(\xi)}$ is the stochastic Jacobian
matrix operator given by the stochastic expansion
\begin{equation}
\mathcal{J(\xi)=}\sum_{\ell=1}^{\widehat{n}}\mathcal{J}_{\ell}\psi_{\ell}%
(\xi). \label{eq:Jacobian-stoch}%
\end{equation}
The expansion is built from matrices $\mathcal{J}_{\ell}\in%
\mathbb{R}
^{n_{x}\times n_{x}}$, $\ell=1,\dots,\widehat{n}$, such that
\[
\mathcal{J}_{1}=\left[
\begin{array}
[c]{cc}%
F_{1} & B^{T}\\
B & 0
\end{array}
\right]  ,\qquad\mathcal{J}_{\ell}=\left[
\begin{array}
[c]{cc}%
F_{\ell} & 0\\
0 & 0
\end{array}
\right]  ,\qquad\ell=2,\dots,\widehat{n},
\]
where $\widehat{n}=\max\left(  n_{\nu},n_{\xi}\right)  $, and $F_{\ell}$ is
a sum of the vector-Laplacian matrix$~A_{\ell}$, the vector-convection
matrix$~N_{\ell}$, and the Newton derivative matrix~$W_{\ell}$,
\[
F_{\ell}=A_{\ell}+N_{\ell}+W_{\ell},
\]
where
\begin{align*}
A_{\ell}  &  =\left[  a_{\ell,ab}\right]  ,\qquad a_{\ell,ab}=\int_{D}%
\nu_{\ell}(x)\,\nabla\phi_{b}:\nabla\phi_{a},\qquad\ell=1,\dots,n_{\nu},\\
N_{\ell}  &  =\left[  n_{\ell,ab}\right]  ,\qquad n_{\ell,ab}=\int_{D}\left(
\vec{u}_{\ell}\cdot\nabla\phi_{b}\right)  \cdot\phi_{a},\qquad\ell
=1,\dots,n_{\xi},\\
W_{\ell}  &  =\left[  w_{\ell,ab}\right]  ,\qquad w_{\ell,ab}=\int_{D}\left(
\phi_{b}\cdot\nabla\vec{u}_{\ell}\right)  \cdot\phi_{a},\qquad\ell
=1,\dots,n_{\xi},
\end{align*}
and if $n_{\nu}>n_{\xi}$, we set $N_{\ell}=W_{\ell}=0$ for $\ell=n_{\xi
+1},\dots,n_{\nu}$, and if $n_{\nu}<n_{\xi}$, we set $A_{\ell}=0$ for
$\ell=n_{\nu+1},\dots,n_{\xi}$. 
\textcolor{black}{
In the numerical experiments, 
we use the stochastic Galerkin method from~\cite{Sousedik-2016-SGM} to calculate the terms~$\vec{u}_\ell$ for the construction of the matrices $N_\ell$.}
The divergence matrix~$B$ is defined as
\[
B=\left[  b_{cd}\right]  ,\qquad b_{cd}=-\int_{D}\varphi_{c}\left(
\nabla\cdot\phi_{d}\right)  ,
\]
and the velocity mass matrix~$G$ is defined as
\[
G=\left[  g_{ab}\right]  ,\qquad g_{ab}=\int_{D}\phi_{b}\,\phi_{a}.
\]

\section{Numerical experiments}

\label{sec:numerical}
\begin{figure}[ptbh]
\centering
\begin{tabular}
[c]{cc}%
\includegraphics[width=6.2cm]{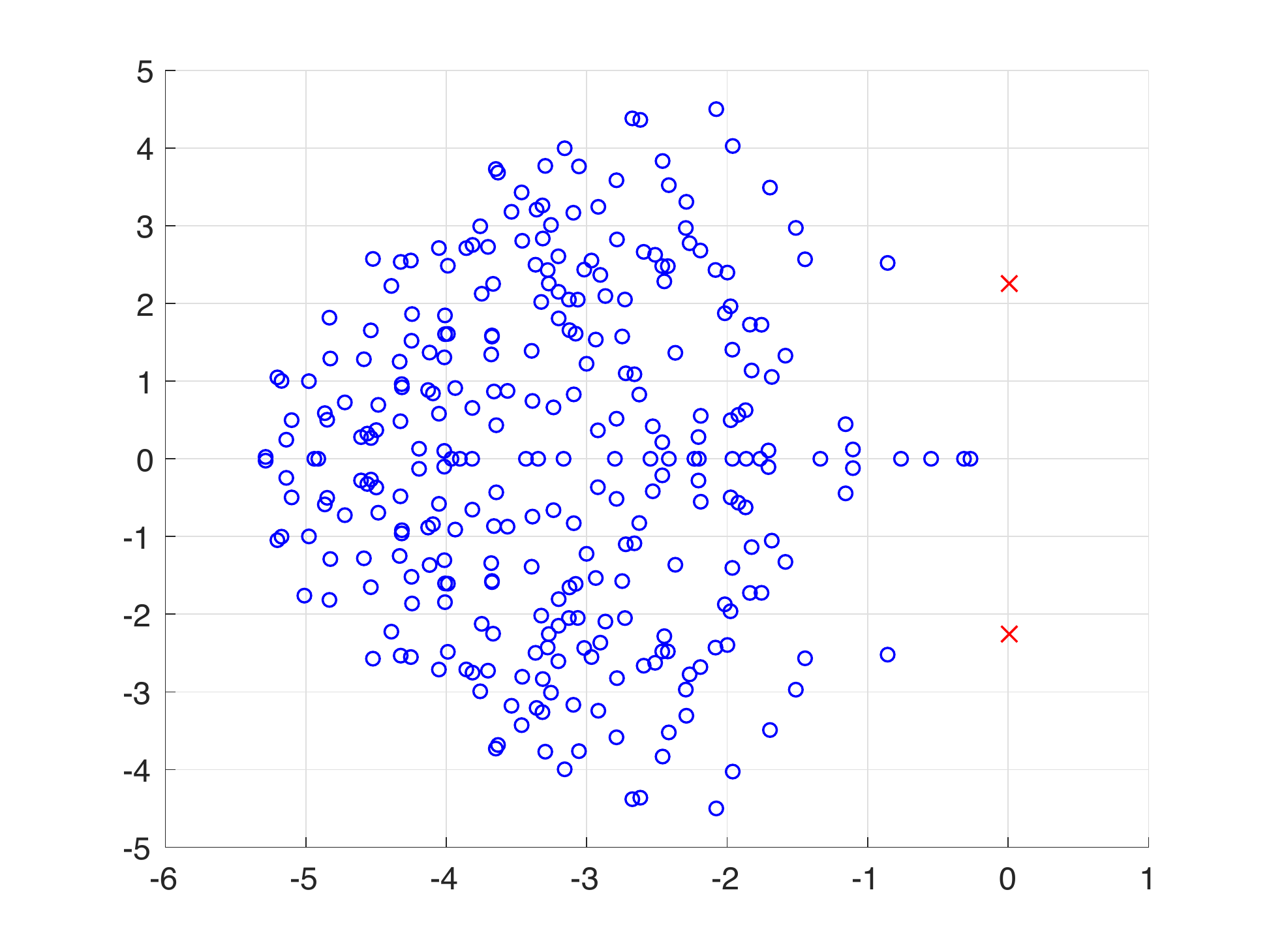} &
\includegraphics[width=6.2cm]{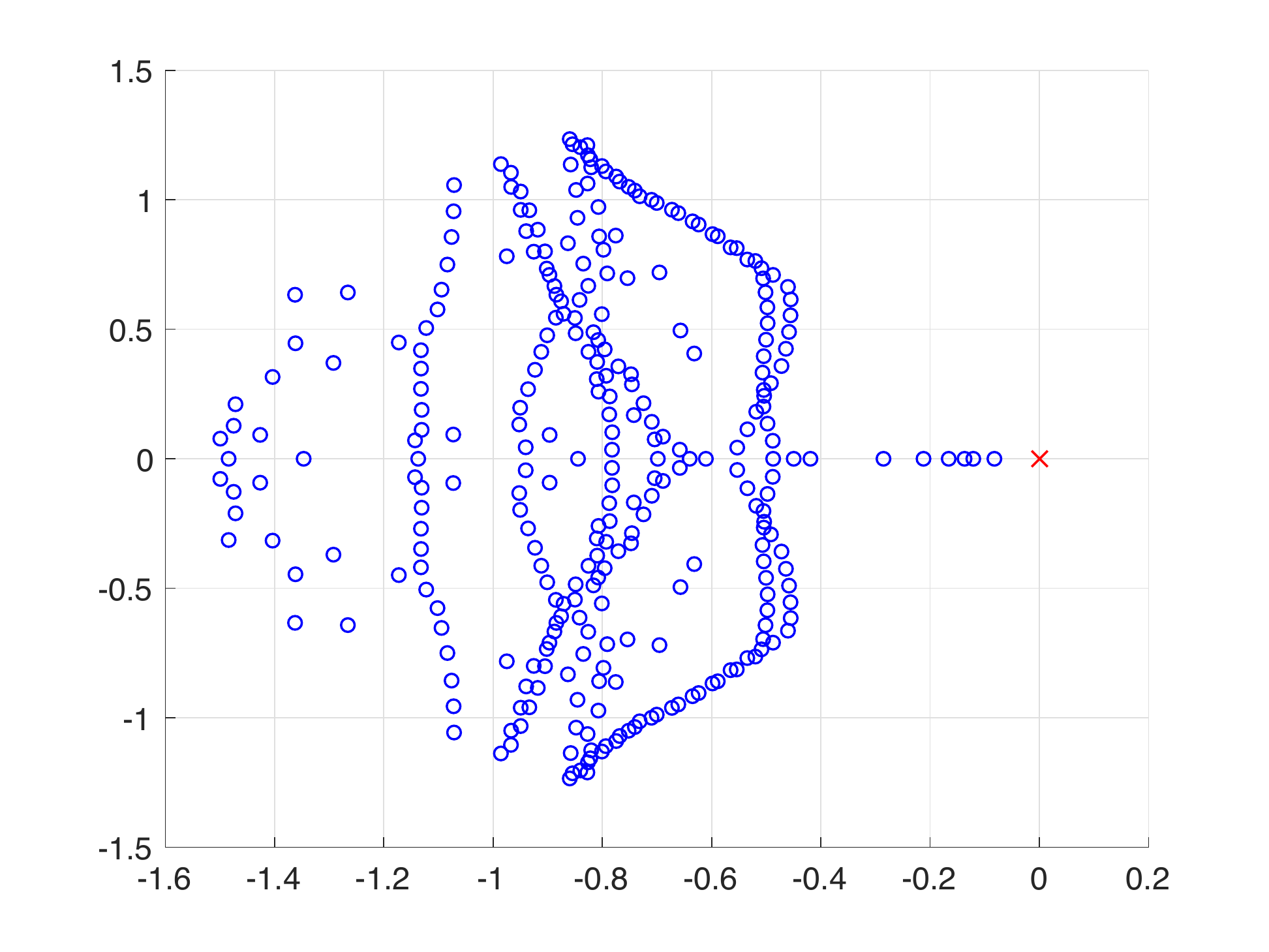}
\end{tabular}
\caption{The $300$ eigenvalues with the largest real part of the
matrices~$K_{1}$ for the two examples: flow around an obstacle (left) and
expansion flow around a symmetric step (right). The rightmost eigenvalues are
indicated by a red cross. }%
\label{fig:lambdaM}%
\end{figure}

We implemented the method in~\textsc{Matlab} using
\textsc{IFISS~3.5} package~\cite{ERS-SIREV}, and we tested the algorithms
using two benchmark problems: flow around an obstacle and an expansion flow
around a symmetric step. The stochastic Galerkin methods were used to solve
both the Navier--Stokes problem (see~\cite{Lee-2019-LRS,Sousedik-2016-SGM} for
full description) and the eigenvalue problem~(\ref{eq:param_eig}), which was
solved using the inexact Newton iteration from Section~\ref{sec:SG}. The
sampling methods (Monte Carlo and stochastic collocation) entail generating a
set of sample viscosities from~(\ref{eq:viscosity}), and for each sample
solving a deterministic steady-state Navier--Stokes equation followed by a
solution of a deterministic eigenvalue problem~(\ref{eq:param_eig_sampling})
with a matrix operator corresponding to sampled Jacobian matrix
operator~(\ref{eq:Jacobian-stoch}), where
the deterministic eigenvalue problems at sample points were solved using
function \texttt{eigs} in \textsc{Matlab}. For the solution of the
Navier--Stokes equation, in both sampling and stochastic Galerkin methods, we
used a hybrid strategy in which an initial approximation was obtained from
solution of the stochastic Stokes problem, after which several steps of Picard
iteration were used to improve the solution, followed by Newton iteration. A
convergent iteration stopped when the Euclidean norm of the algebraic residual
was smaller than~$10^{-8}$, see~\cite{Sousedik-2016-SGM} for more details.
Also, when replacing the mass matrix$~M$ by the shifted mass matrix$~M_{\sigma
}$, see (\ref{eq:matrices-det}) and (\ref{eq:M-shifted}), we set
$\sigma=-10^{-2}$ as in~\cite{Elman-2012-LII}. The $300$ eigenvalues with the
largest real part of the deterministic eigenvalue problem with mean
viscosity~$\nu_{1}$ for each of the two examples are displayed in
Figure~\ref{fig:lambdaM}.

\subsection{Flow around an obstacle}

\label{sec:obstacle}
\begin{figure}[ptbh]
\begin{center}
\includegraphics[width=9.8cm]{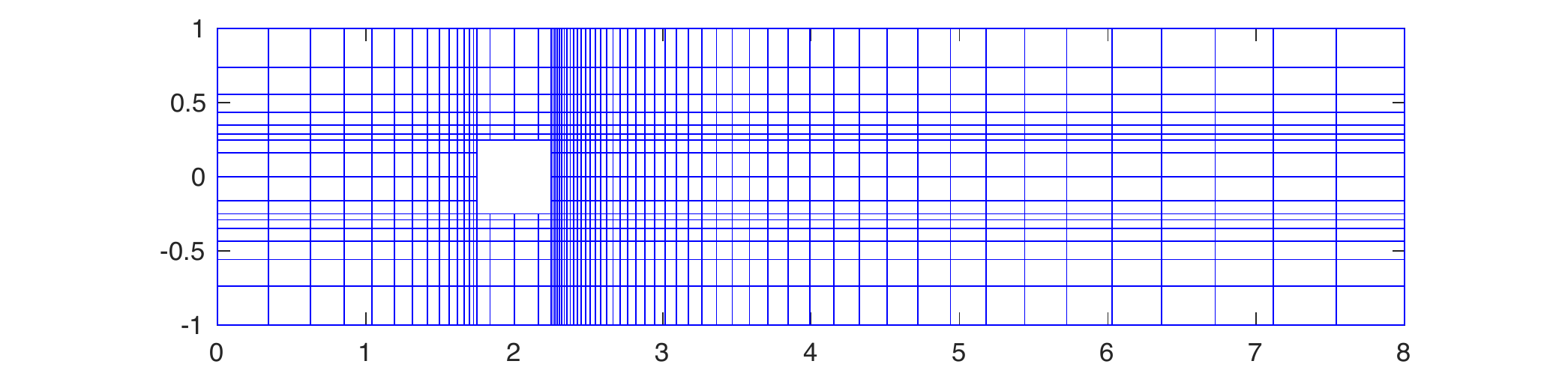}
\end{center}
\caption{Finite element mesh for the flow around an obstacle problem. }%
\label{fig:mesh-obstacle}%
\end{figure}

\begin{figure}[ptbh]
\centering
\begin{tabular}
[c]{cc}%
\includegraphics[width=6.3cm]{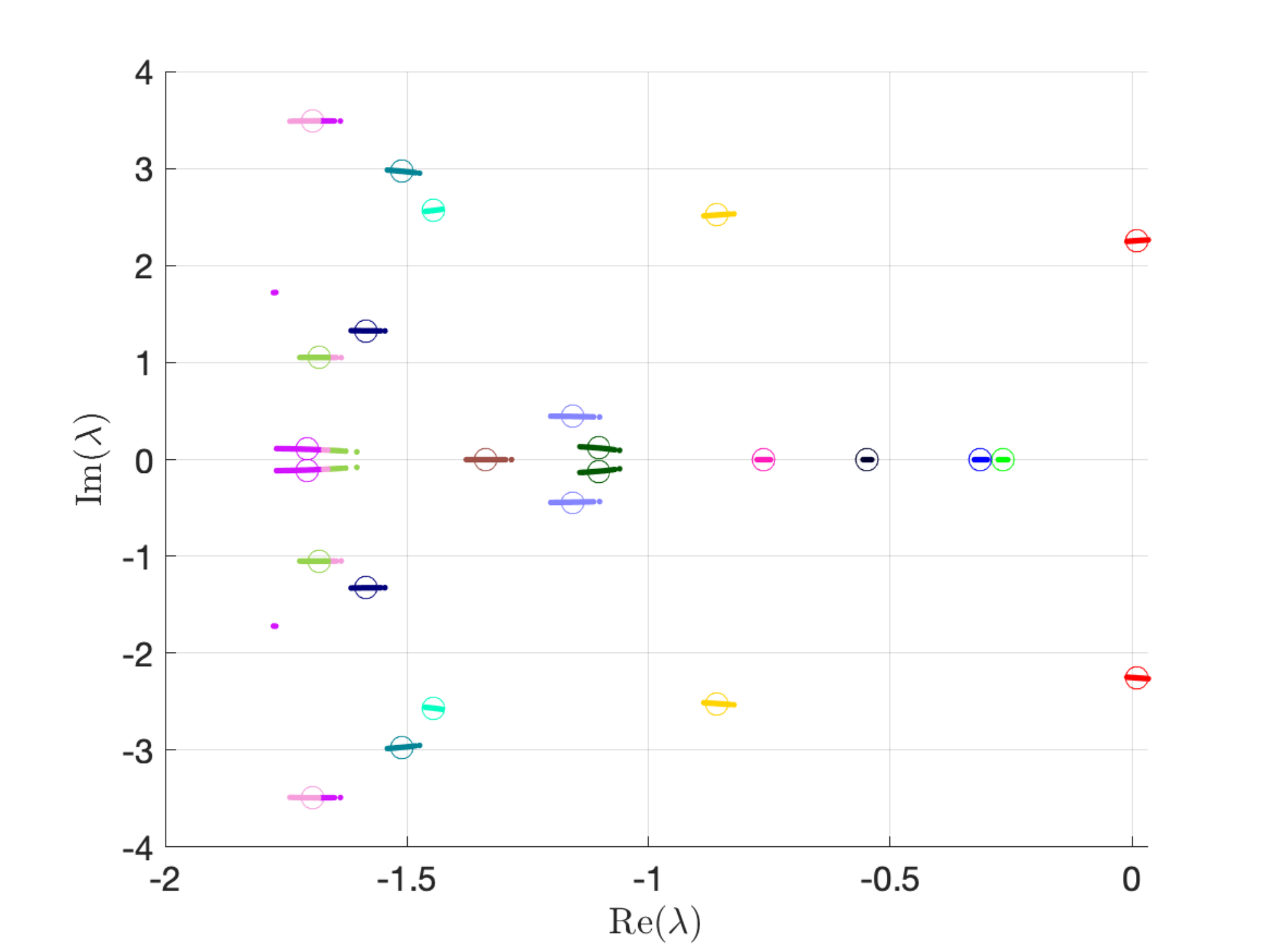} &
\includegraphics[width=6.3cm]{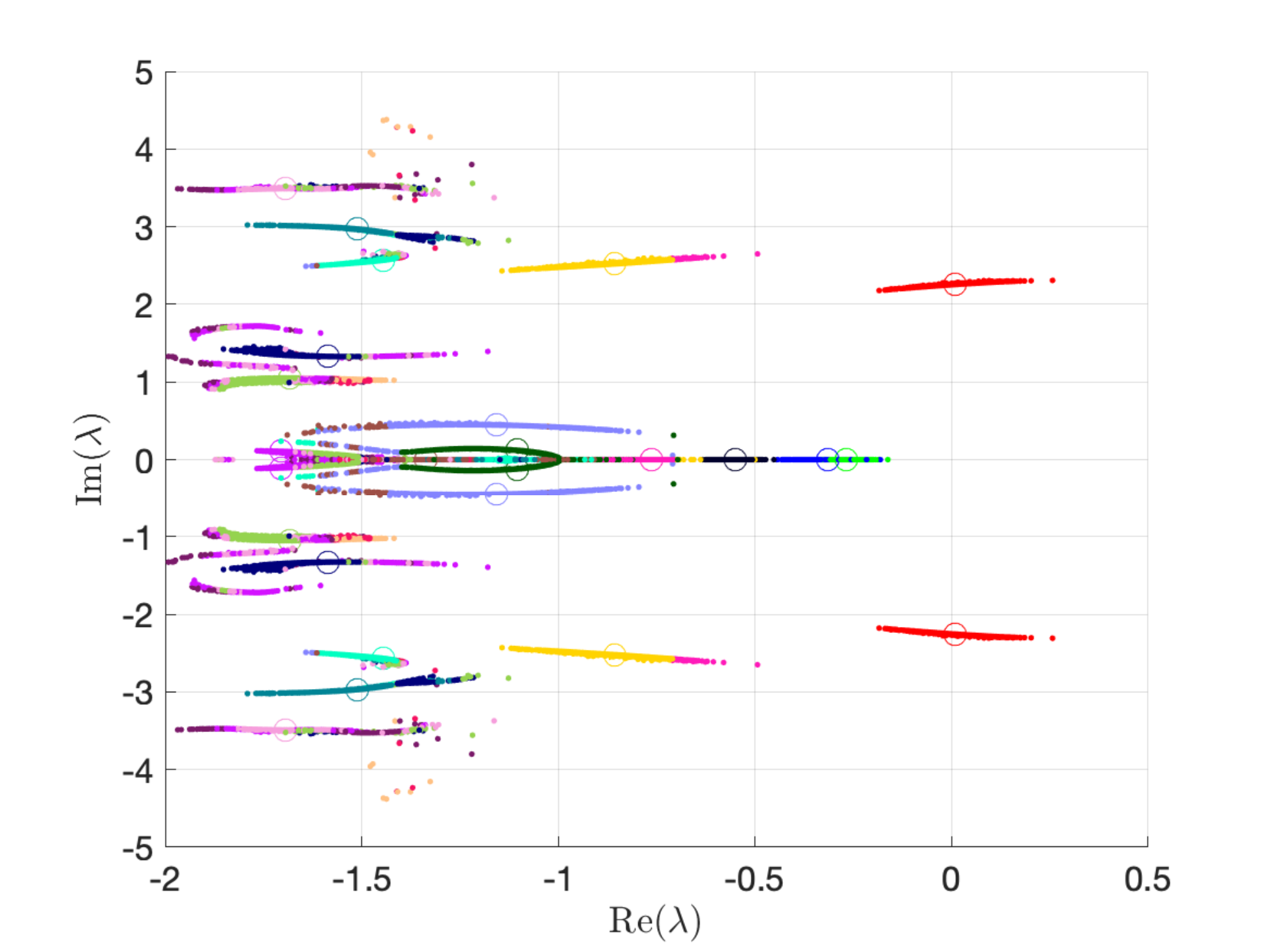}
\end{tabular}
\caption{Monte Carlo samples of~$25$ eigenvalues with the largest real part
for the flow around an obstacle with $CoV=1\%$ (left) and $CoV=10\%$ (right).
The eigenvalues of the mean problem are indicated by circles. }%
\label{fig:obstacleMC}%
\end{figure}

\begin{figure}[ptbh]
\centering
\begin{tabular}
[c]{cc}%
\includegraphics[width=6.5cm]{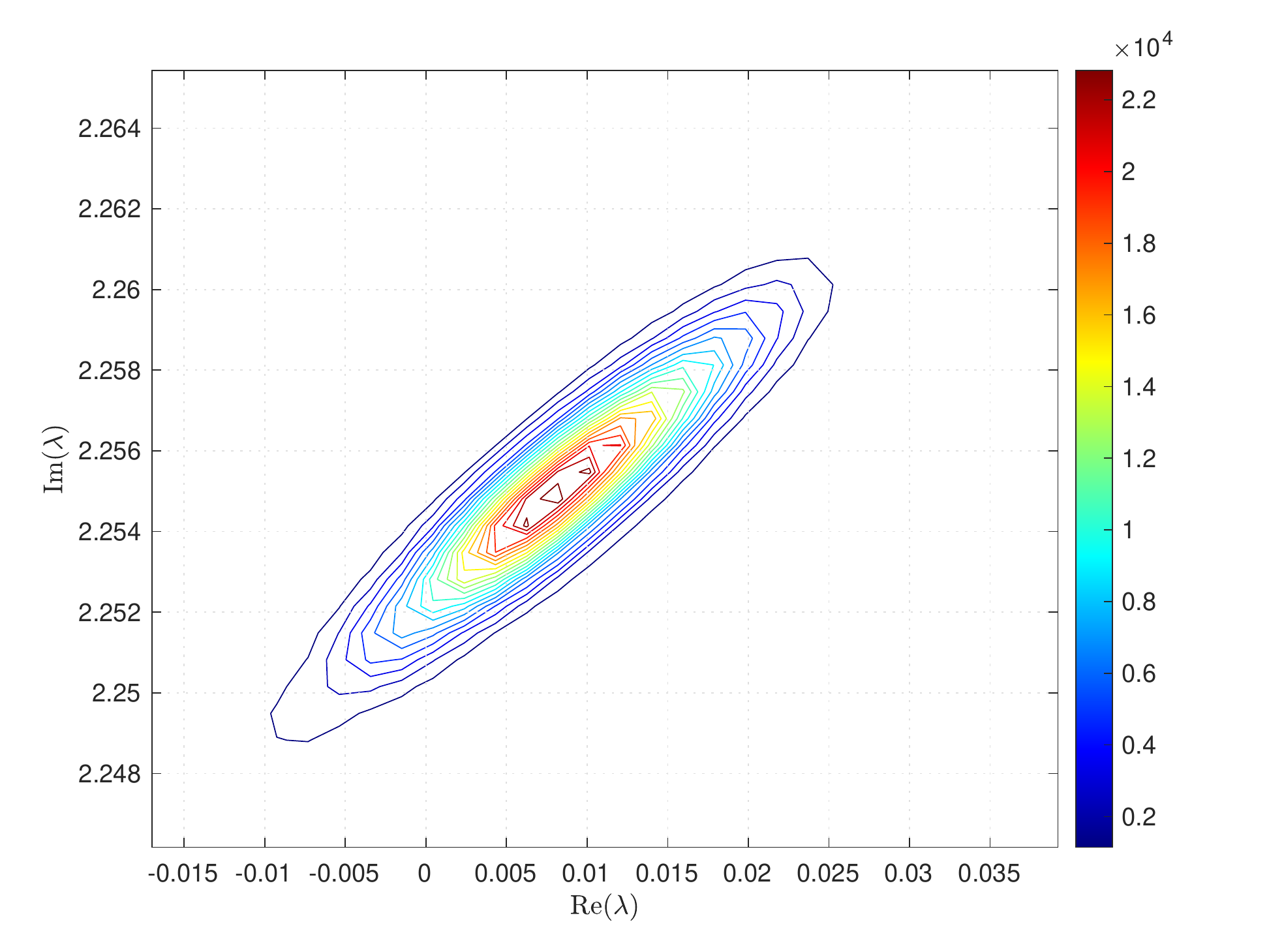} &
\includegraphics[width=6.5cm]{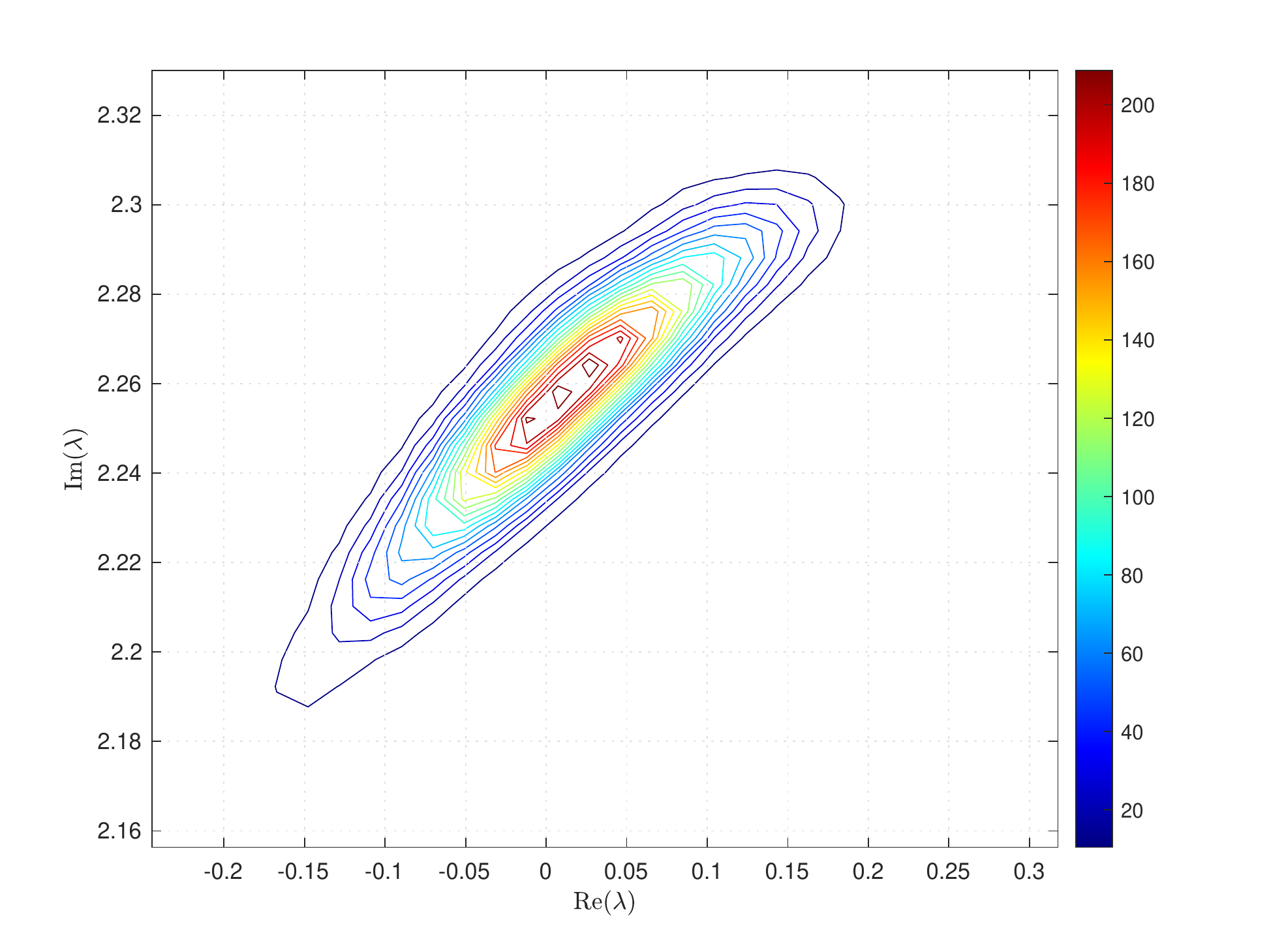}\\
\includegraphics[width=6.5cm]{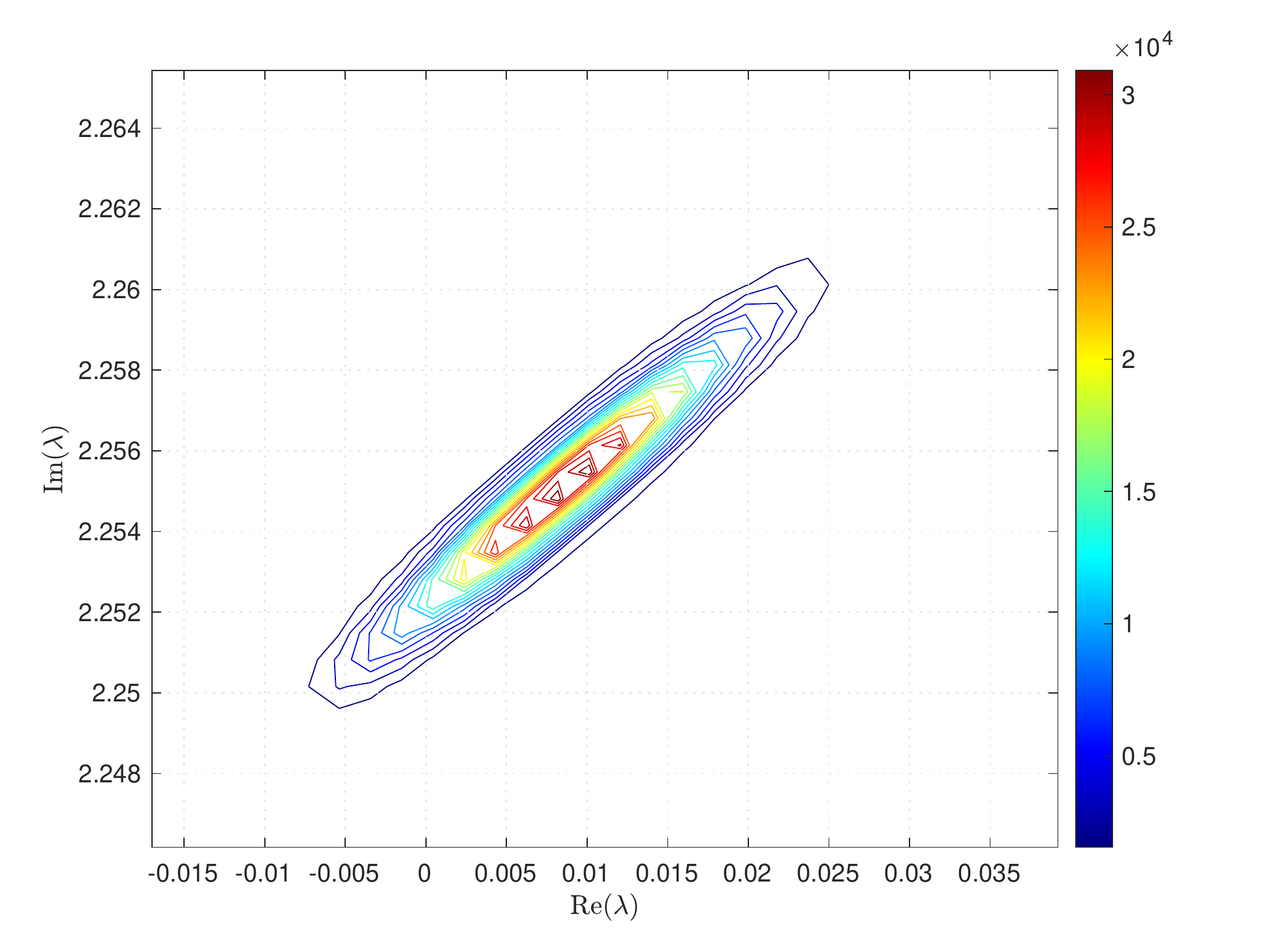} &
\includegraphics[width=6.5cm]{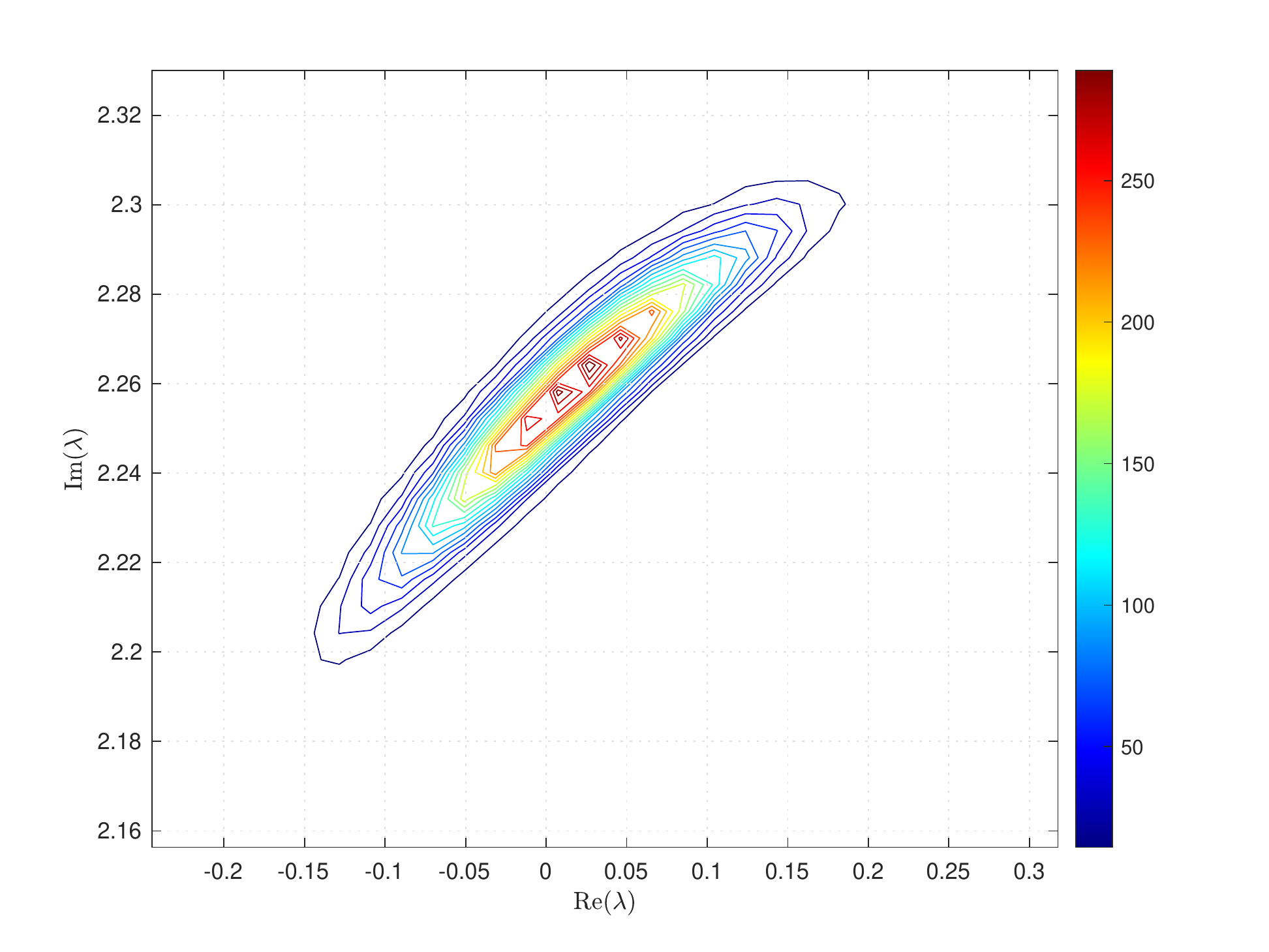}\\
\includegraphics[width=6.5cm]{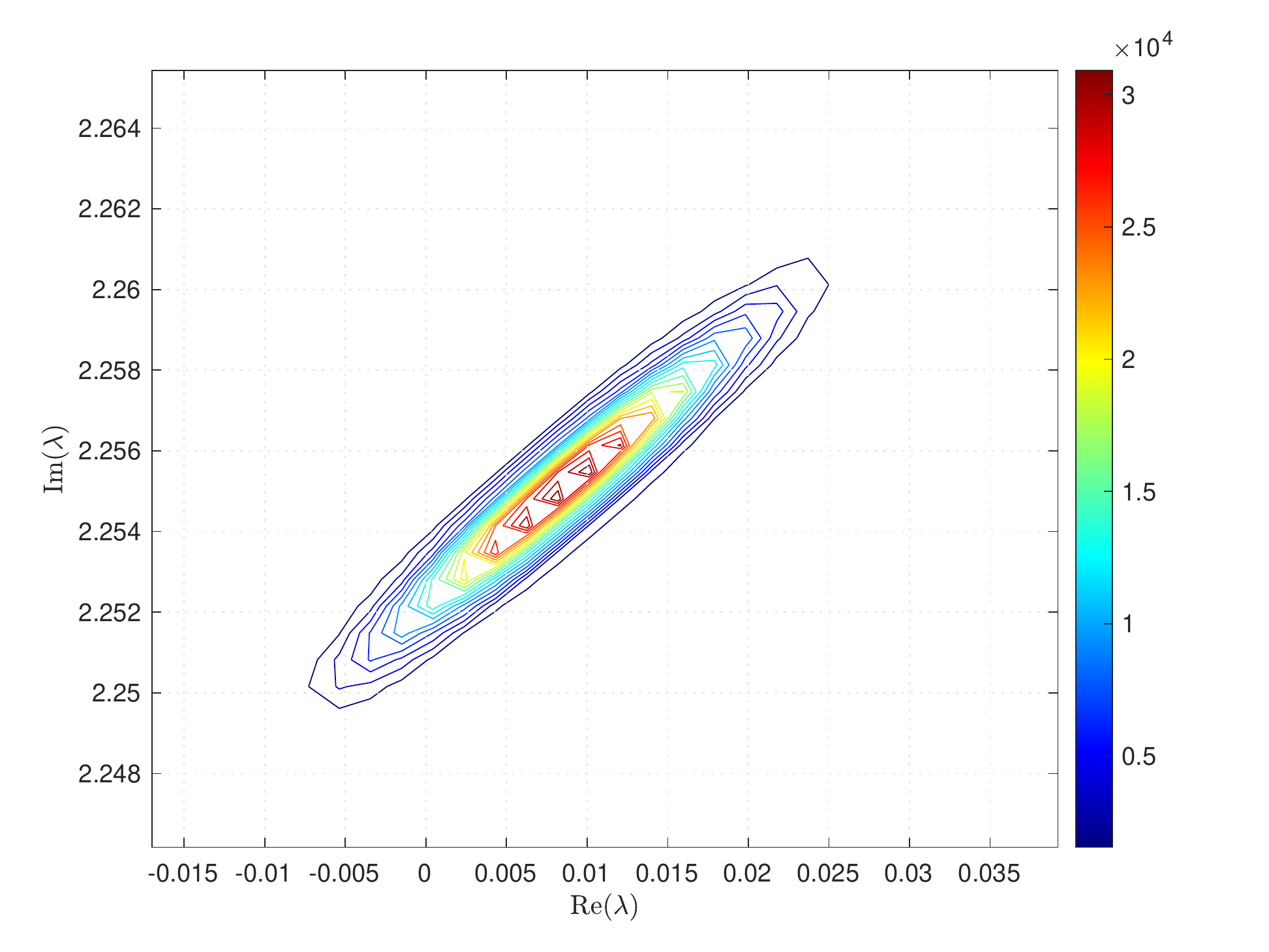} &
\includegraphics[width=6.5cm]{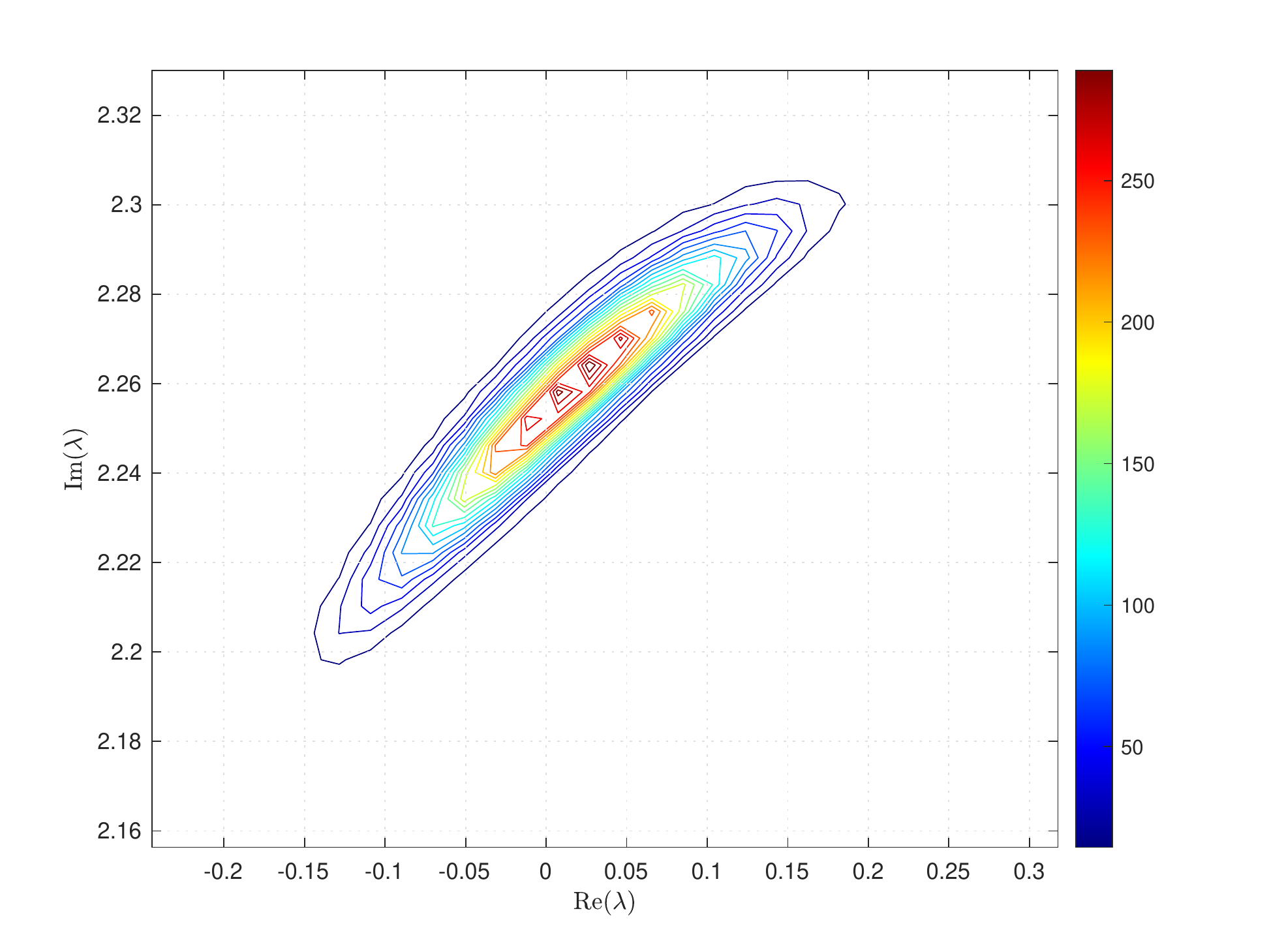}
\end{tabular}
\caption{Plots of the pdf estimate of the rightmost eigenvalue with positive
imaginary part obtained using Monte Carlo (top), stochastic collocation
(middle) and stochastic Galerkin method (bottom) for the flow around an
obstacle with $CoV=1\%$ (left) and $CoV=10\%$ (right).}%
\label{fig:obstacle-pdf2D}%
\end{figure}

\begin{figure}[ptbh]
\centering
\begin{tabular}
[c]{cc}%
\includegraphics[width=6.3cm]{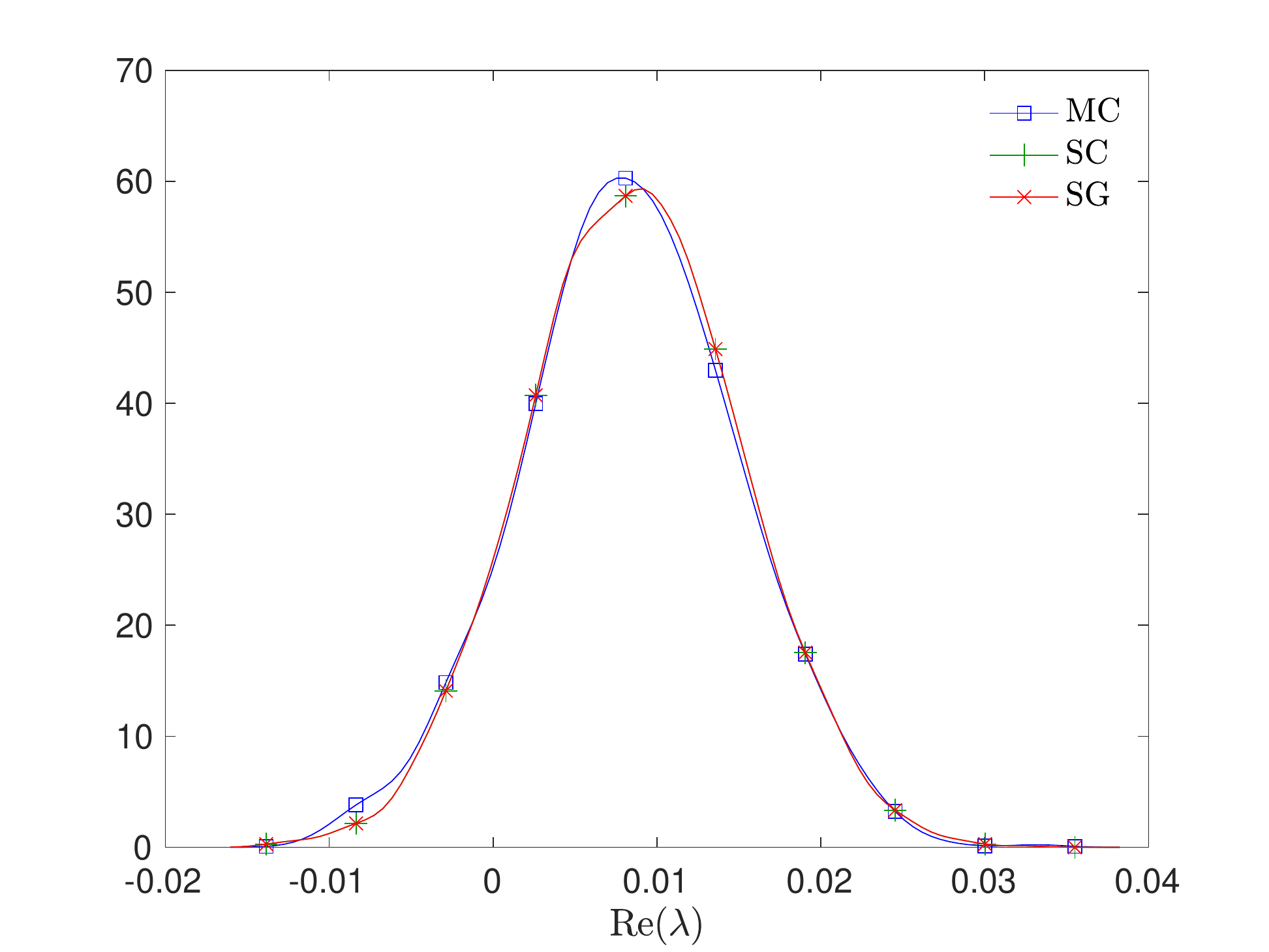} &
\includegraphics[width=6.3cm]{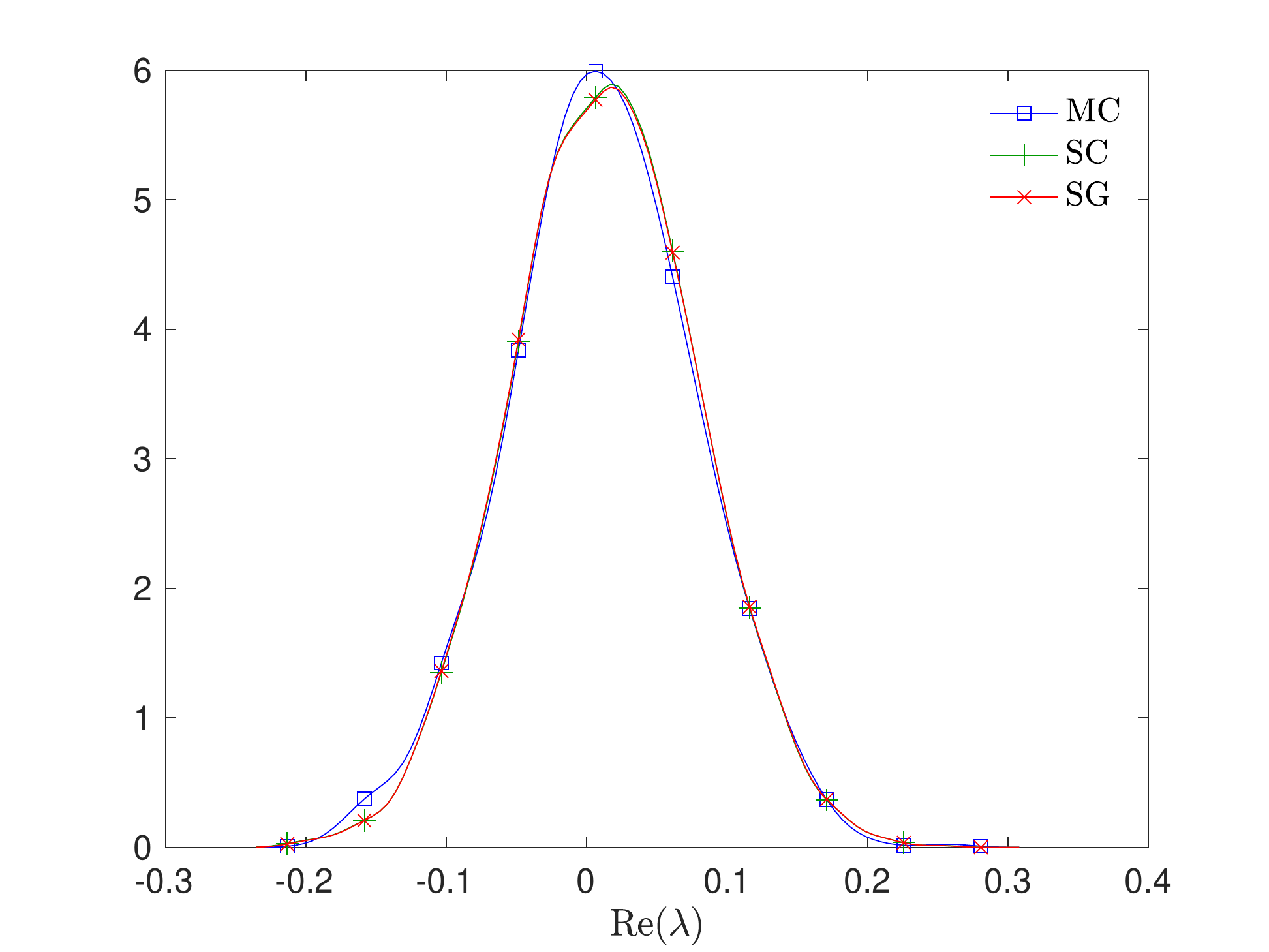}
\end{tabular}
\caption{Plots of the pdf estimate of the real part of the rightmost
eigenvalue obtained using Monte Carlo (MC), stochastic collocation (SC) and
stochastic Galerkin method (SG) for the flow around an obstacle with $CoV=1\%$
(left) and $CoV=10\%$ (right).}%
\label{fig:obstacle-pdf-Re}%
\end{figure}

For the first example, we considered flow around an
obstacle in a similar setup as studied in~\cite{Sousedik-2016-SGM}.
The domain of the channel and the discretization are shown in
Figure~\ref{fig:mesh-obstacle}. The spatial discretization used a stretched
grid with $1008$ \textbf{\textit{Q}}$_{2}-$\textbf{\textit{Q}}$_{1}$
finite\ elements. We note that these elements are referred to as
\emph{Taylor-Hood} in the literature.
There were $8416$ velocity and $1096$ pressure degrees of freedom. The
viscosity$~\nu(x,\xi)$\ was taken to be a truncated lognormal process
transformed from the underlying Gaussian process~\cite{Ghanem-1999-NGS}. That
is, $\psi_{\ell}(\xi)$, $\ell=1,\dots,n_{\nu}$, is a set of Hermite
polynomials and, denoting the coefficients of the Karhunen-Lo\`{e}ve expansion
of the Gaussian process by$~g_{j}(x)$ and $\eta_{j}=\xi_{j}-g_{j},$
$j=1,\dots,m_{\xi}$, the coefficients in expansion~(\ref{eq:viscosity}) were
computed as
\[
\nu_{\ell}(x)=\frac{\mathbb{E}\left[  \psi_{\ell}(\eta)\right]  }%
{\mathbb{E}\left[  \psi_{\ell}^{2}(\eta)\right]  }\exp\left[  g_{0}+\frac
{1}{2}\sum_{j=1}^{m_{\xi}}\left(  g_{j}(x)\right)  ^{2}\right]  .
\]
The covariance function of the Gaussian field, for points $X_{1}=(x_{1}%
,y_{1})$ and $X_{2}=(x_{2},y_{2})$ in$~D$, was chosen to be
\begin{equation}
C\left(  X_{1},X_{2}\right)  =\sigma_{g}^{2}\exp\left(  -\frac{\left\vert
x_{2}-x_{1}\right\vert }{L_{x}}-\frac{\left\vert y_{2}-y_{1}\right\vert
}{L_{y}}\right)  ,\label{eq:covariance_kernel}%
\end{equation}
where$~L_{x}$ and $L_{y}$ are the correlation lengths of the random variables
$\xi_{i}$, $i=1,\dots,m_{\xi}$, in the $x$ and $y$\ directions, respectively,
and $\sigma_{g}$ is the standard deviation of the Gaussian random field. The
correlation lengths were set to be equal to $25\%$ of the width and height of
the domain. The coefficient of variation $CoV$ of the lognormal field, defined
as $CoV=\sigma_{\nu}/\nu_{1}$, where $\sigma_{\nu}$ is the standard deviation
and $\nu_{1}$ is the mean viscosity, was $1\%$ or $10\%$. The
viscosity~(\ref{eq:viscosity}) was parameterized using $m_{\xi}=2$ random
variables. According to~\cite{Matthies-2005-GML}, in order to guarantee a
complete representation of the lognormal process by~(\ref{eq:viscosity}), the
degree of polynomial expansion of~$\nu(x,\xi)$ should be twice the degree of
the expansion of the solution. We followed the same strategy here. Therefore,
the values of $n_{\xi}$ and $n_{\nu}$ are, cf., e.g.~\cite[p.~87]%
{Ghanem-1991-SFE} or~\cite[Section~5.2]{Xiu-2010-NMS}, $n_{\xi}=\frac{\left(
m_{\xi}+p\right)  !}{m_{\xi}!p!}$, $n_{\nu}=\frac{\left(  m_{\xi}+2p\right)
!}{m_{\xi}!\left(  2p\right)  !}$. For the gPC expansion of
eigenvalues/eigenvectors~(\ref{eq:sol_mat}), the maximal degree of
gPC\ expansion is $p=3$, so then $n_{\xi}=10$ and $n_{\nu}=28$. We used
$1\times10^{3}$ samples for the Monte Carlo method
and Smolyak sparse grid with Gauss-Hermite quadrature points and grid
level$~4$ for collocation
\textcolor{black}{see, e.g.,~\cite{LeMaitre-2010-SMU}\ for discussion of quadrature rules}.
With these settings, the size of $\left\{  H_{\ell}\right\}  _{\ell=1}%
^{n_{\nu}}$\ in~(\ref{eq:H}) was $10\times10\times28$\ with $203$ nonzeros,
and there were $n_{q}=29$ points in the sparse grid.
\textcolor{black}{As a consequence, the size of the stochastic Galerkin matrices is $n_\xi (n_u+n_p)=95120$,
the matrix associated with the Jacobian is fully block dense and the mass matrix is block diagonal,
but we note that these matrices are never formed in implementation.} For the
solution of the Navier--Stokes problem we used the hybrid strategy with $6$
steps of Picard iteration followed by at most $15$ steps of Newton iteration.
We used mean viscosity $\nu_{1}=5.36193\times10^{-3}$, which corresponds to
Reynolds number $Re=373$, and the rightmost eigenvalue pair is $0.0085\pm
2.2551i$, see the left panel in Figure~\ref{fig:lambdaM}. The
Figure~\ref{fig:obstacleMC} displays Monte Carlo realizations of the $25$
eigenvalues with the largest real part for the values $CoV=1\%$ and
$CoV=10\%$. It can be seen that the rightmost eigenvalue is relatively less
sensitive to perturbation comparing to the other eigenvalues, and since its
real part is well separated from the rest of the spectrum, it can be easily
identified in all runs of a sampling method.
Figure~\ref{fig:obstacle-pdf2D} displays the probability density function
(pdf) estimates of the rightmost eigenvalue with the positive imaginary part
obtained directly by Monte Carlo, the stochastic collocation and stochastic
Galerkin methods, for which the estimates were obtained using \textsc{Matlab}
function \texttt{ksdensity} (in 2D) for sampled gPC\ expansions.
Figure~\ref{fig:obstacle-pdf-Re} shows plots of the estimated pdf of the real
part of the rightmost eigenvalue. In both figures we can see a good agreement
of the plots in the left column corresponding to $CoV=1\%$ and in the right
column corresponding to $CoV=10\%$.\ In Table~\ref{tab:obstacle_gPC} we
tabulate the coefficients of the gPC\ expansion of the rightmost eigenvalue
with positive imaginary part computed using the stochastic collocation and
Galerkin methods. A good agreement of coefficients can be seen, in particular
for coefficients with value much larger than zero, specifically with
$k=1,2,4,6,7$ and$~9$. Finally, in Table~\ref{tab:gmres_obstacle} we examine
the inexact line-search Newton iteration from Algorithm~\ref{alg:line_search}.
For the line-search method, we set $\rho=0.9$ for the backtracking and
$c=0.25$. The initial guess is set using the rightmost eigenvalue and
corresponding eigenvector of the eigenvalue
problem~(\ref{eq:eigenproblem-mean}) concatenated by zeros. The nonlinear
iteration terminates when the norm of the residual $\left\Vert \widehat{r}%
_{n}\right\Vert _{2}<10^{-10}$.\ The linear systems in Line~\ref{ln:lin_sys}
of Algorithm~\ref{alg:line_search} are solved using GMRES\ with the mean-based
preconditioner (Algorithm~\ref{alg:MB}), constraint mean-based preconditioner
(Algorithm~\ref{alg:cMB}) and its updated variant discussed in
Section~\ref{sec:SMW}, and the constraint hierarchical Gauss-Seidel
preconditioner (Algorithm~\ref{alg:chGS}--\ref{alg:chGS_cont}), which was used
without truncation of the matrix-vector multiplications and also with
truncation, setting $p_{t}=2$, as discussed in
Section~\ref{sec:preconditioners}. For the mean-based preconditioner we used
$\epsilon_{\operatorname{Re}}=\epsilon_{\operatorname{Im}}=0.97$, which worked
best in our experience, and $\epsilon_{\operatorname{Re}}=\epsilon
_{\operatorname{Im}}=1$\ otherwise. For the constraint mean-based
preconditioner the matrix$~\mathcal{M}_{\text{cMB}}$
from~(\ref{eq:M-preconditioner}) was factored using LU decomposition, and the
updated variant from Section~\ref{sec:SMW} was used also in the constraint
hierarchical Gauss-Seidel preconditioner. First, we note that the performance
of the algorithm with the mean-based preconditioner was very sensitive to the
choice of $\epsilon_{\operatorname{Re}}$ and $\epsilon_{\operatorname{Im}}$,
and it can be seen that it is quite sensitive also to $CoV$. On the other
hand, the performance of all variants of the constraint preconditioners appear
to be far less sensitive, and we see only a mild increase in numbers of both
nonlinear and GMRES\ iterations. Next, we see that updating the constraint
mean-based preconditioner reduces the numbers of GMRES iterations in
particular in the latter\ steps of the nonlinear method. Finally, we see that
using the constraint hierarchical Gauss-Seidel preconditioner further
decreases the number of GMRES\ iterations, for smaller $CoV$ it may be
suitable to truncate the matrix-vector multiplications without any change in
iteration counts and even though we see with $CoV=10\%$ an increase in number
of nonlinear steps, the overall number of GMRES\ iterations remains smaller
than when the two variants of the constraint mean-based preconditioner were used.

\begin{table}[b]
\caption{The $10$ coefficients of the gPC expansion of the rightmost
eigenvalue with positive complex part for the flow around an obstacle problem
with $CoV=1\%$ and $10\%$ computed using stochastic collocation (SC), and
stochastic Galerkin method (SG). Here $d$ is the polynomial degree and $k$ is
the index of basis function in expansion~(\ref{eq:sol_mat}).}%
\label{tab:obstacle_gPC}
\begin{center}
{\footnotesize \renewcommand{\arraystretch}{1.3}
\begin{tabular}
[c]{|c|c|cc|}\hline
$d$ & $k$ & SC & SG\\\hline
&  & \multicolumn{2}{|c|}{$CoV=1\%$}\\\hline
0 & 1 & 8.5726E-03 + 2.2551E+00\,i & 8.5726E-03 + 2.2551E+00\,i\\\hline
\multirow{2}{*}{1} & 2 & -6.5686E-03 - 2.2643E-03\,i & -6.5686E-03 -
2.2643E-03\,i\\
& 3 & 1.1181E-16 - 2.0817E-14\,i & 2.6512E-17 + 8.3094E-17\,i\\\hline
\multirow{3}{*}{2} & 4 & -1.1802E-06 - 2.4274E-05\,i & -1.2055e-06 -
2.4200E-05\,i\\
& 5 & 3.8351E-15 - 4.4964E-15\,i & 8.9732E-20 - 2.0565E-19\,i\\
& 6 & -3.3393E-06 + 4.0603E-05\,i & -3.3527E-06 + 4.0641E-05\,i\\\hline
\multirow{4}{*}{3} & 7 & -1.0635E-07 + 4.1735E-07\,i & -8.5671E-08 +
3.5926E-07\,i\\
& 8 & 7.8095E-16 + 6.1617E-15\,i & -4.3191E-22 - 8.3970E-21\,i\\
& 9 & -4.6791E-07 + 5.1602E-08\,i & -4.4762E-07 - 6.0766E-09\,i\\
& 10 & 2.2155E-15 + 4.6907E-15\,i & 1.2691E-15 + 2.9181E-16\,i\\\hline
&  & \multicolumn{2}{|c|}{$CoV=10\%$}\\\hline
0 & 1 & 1.3420E-02 + 2.2577E+00\,i & 1.3419E-02 + 2.2576E+00\,i\\\hline
\multirow{2}{*}{1} & 2 & -6.6200E-02 - 2.2034E-02\,i & -6.6243E-02 -
2.2018E-02\,i\\
& 3 & 1.6011E-15 - 1.0297E-14\,i & 1.1672E-15 + 8.8396E-16\,i\\\hline
\multirow{3}{*}{2} & 4 & -2.2415E-04 - 2.5416E-03\,i & -1.0889E-04 -
2.4178E-03\,i\\
& 5 & 8.5869E-17 - 1.0547E-15\,i & 1.1865E-17 + 6.5559E-17\,i\\
& 6 & -2.7323E-04 + 4.1219E-03\,i & -2.1977E-04 + 4.1437E-03\,i\\\hline
\multirow{4}{*}{3} & 7 & -4.8106E-05 + 3.556E-04\,i & 1.3510E-04 +
9.1486E-05\,i\\
& 8 & 2.8365E-15 + 6.1062E-15\,i & 8.0683E-19 + 5.3753E-18\,i\\
& 9 & -4.5696E-04 + 2.7795E-06\,i & -4.1149E-04 - 1.8160E-04\,i\\
& 10 & 1.7408E-15 + 1.3101E-14\,i & 1.3975E-15 + 3.5152E-16\,i\\\hline
\end{tabular}
}
\end{center}
\end{table}

\begin{table}[b]
\caption{The number of GMRES iterations in each step of inexact line-search
Newton method (Algorithm~\ref{alg:line_search}) for computing the rightmost
eigenvalue and corresponding eigenvectors of the flow around an obstacle
problem with $CoV = 1\%$ (left) and $10\%$ (right) and with the stopping
criteria $\Vert r_{n}\Vert_{2}<10^{-10}$ and different choices of
preconditioners: mean-based (MB) from Algorithm~\ref{alg:MB}, constraint
mean-based (cMB) from Algorithm~\ref{alg:cMB} and its updated variant
(cMB$_{u}$) from Section~\ref{sec:SMW}, and the constraint hierarchical
Gauss-Seidel preconditioner (chGS) from Algorithm~\ref{alg:chGS}%
--\ref{alg:chGS_cont} and also with truncation, setting $p_{t}=2$ (chGS$_{2}%
$). }%
\label{tab:gmres_obstacle}
\begin{center}
{\footnotesize \renewcommand{\arraystretch}{1.3}
\begin{tabular}
[c]{|c|ccccc|ccccc|}\hline
& \multicolumn{5}{|c|}{$CoV=1\%$} & \multicolumn{5}{c|}{$CoV=10\%$}\\\hline
step & MB & cMB & cMB$_{u}$ & chGS & chGS$_{2}$ & MB & cMB & cMB$_{u}$ &
chGS & chGS$_{2}$\\\hline
1 & 2 & 1 & 1 & 1 & 1 & 7 & 1 & 1 & 1 & 1\\
2 & 2 & 1 & 1 & 1 & 1 & 6 & 3 & 2 & 3 & 3\\
3 & 6 & 3 & 2 & 1 & 1 & 13 & 4 & 4 & 3 & 4\\
4 & 9 & 6 & 3 & 2 & 2 & 10 & 8 & 7 & 3 & 4\\
5 & 15 & 10 & 6 & 3 & 3 & 15 & 16 & 13 & 4 & 5\\
6 &  &  &  &  &  & 14 &  &  & 8 & 8\\
7 &  &  &  &  &  & 25 &  &  &  & \\
8 &  &  &  &  &  & 32 &  &  &  & \\
9 &  &  &  &  &  & 67 &  &  &  & \\\hline
\end{tabular}
}
\end{center}
\end{table}

\subsection{Expansion flow around a symmetric step}

\label{sec:expansion}
\begin{figure}[ptbh]
\begin{center}
\includegraphics[width=13cm]{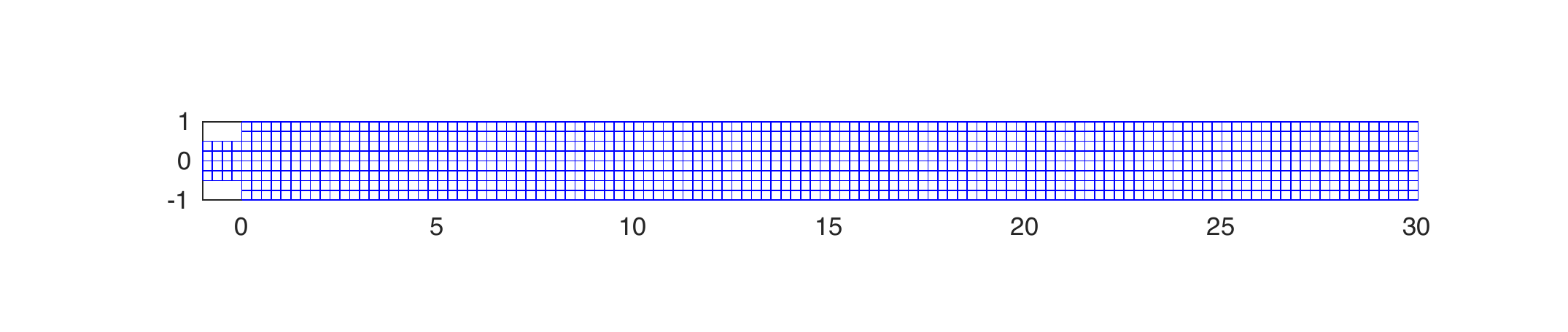}
\end{center}
\caption{Finite element mesh for the expansion flow around a symmetric step. }%
\label{fig:mesh-symstep}%
\end{figure}

\begin{figure}[ptbh]
\centering
\begin{tabular}
[c]{cc}%
\includegraphics[width=6.3cm]{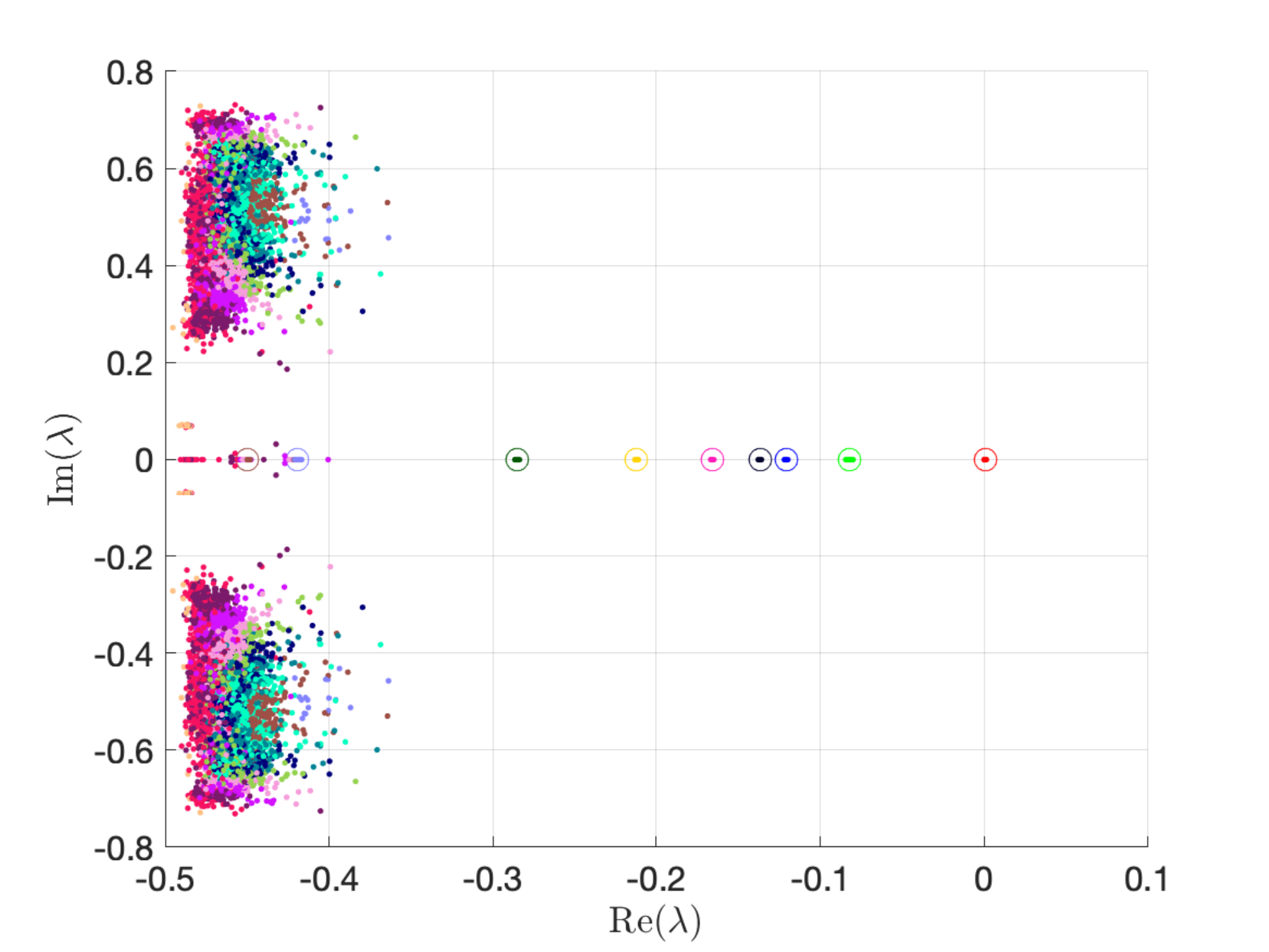} &
\includegraphics[width=6.3cm]{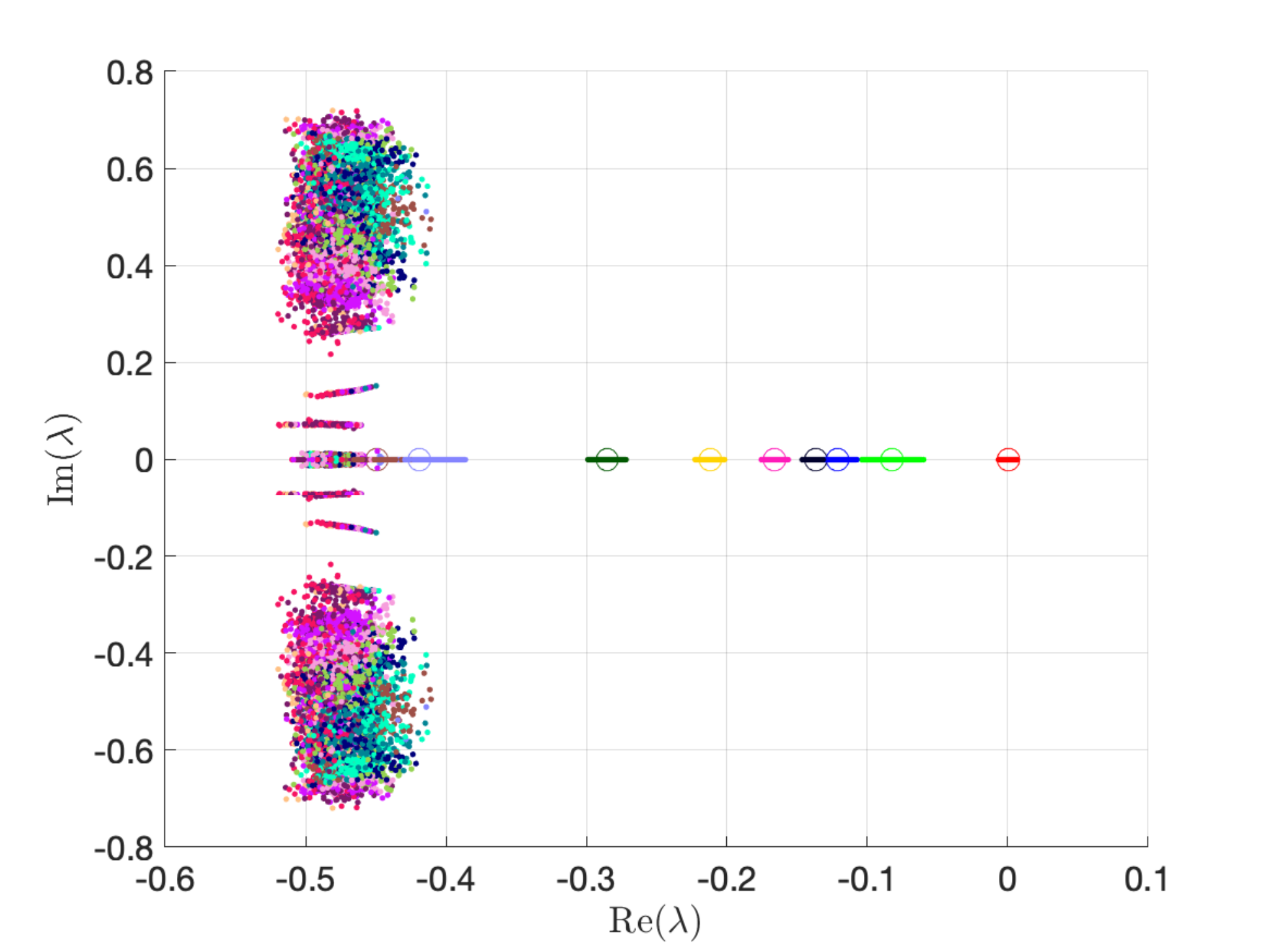}
\end{tabular}
\caption{Monte Carlo samples of~$25$ eigenvalues with the largest real part
for the flow around an obstacle with $CoV=1\%$ (left) and $CoV=10\%$ (right).
The eigenvalues of the mean problem are indicated by circles. }%
\label{fig:symstepMC}%
\end{figure}

\begin{figure}[ptbh]
\centering
\begin{tabular}
[c]{cc}%
\includegraphics[width=6.3cm]{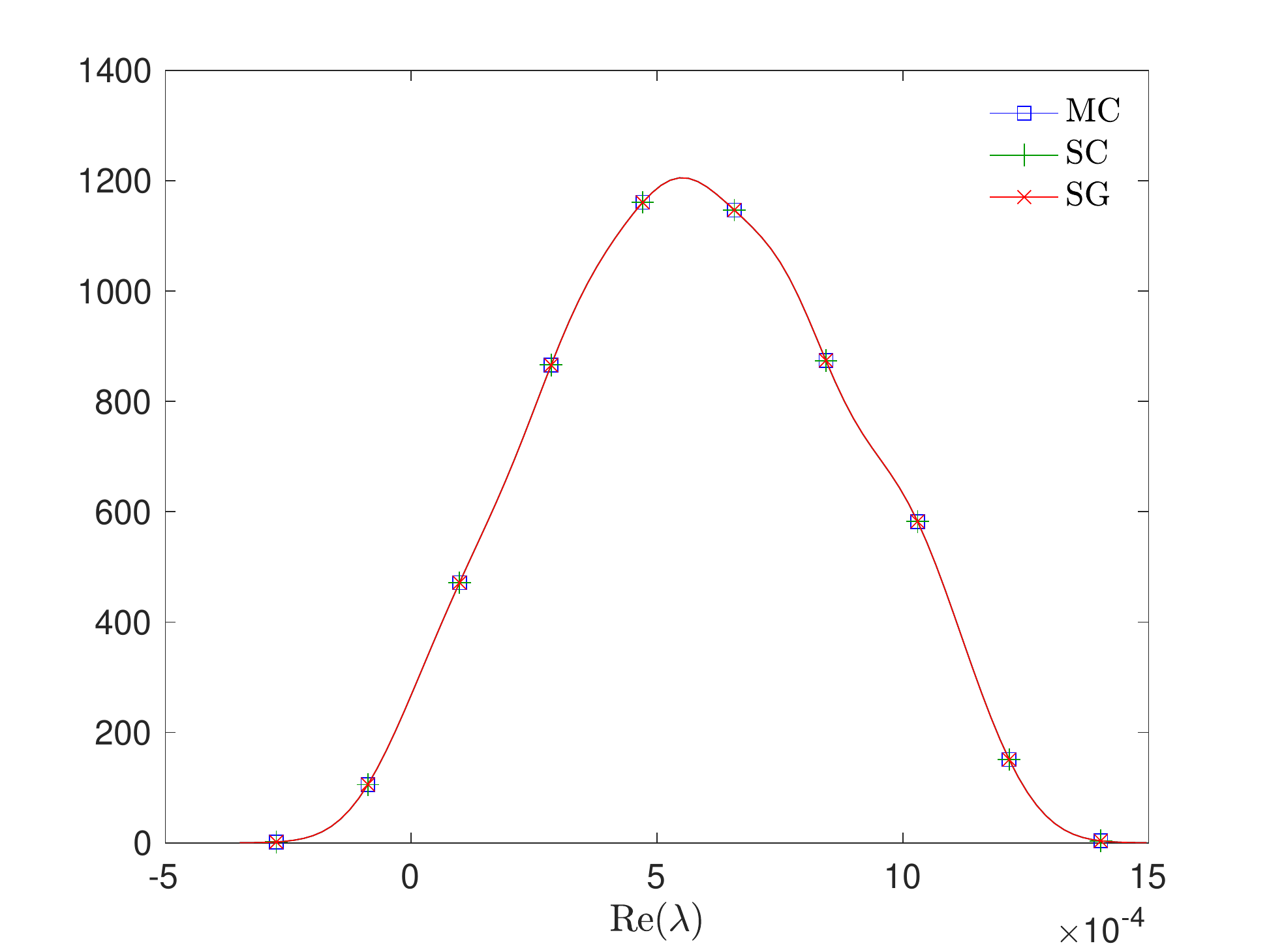} &
\includegraphics[width=6.3cm]{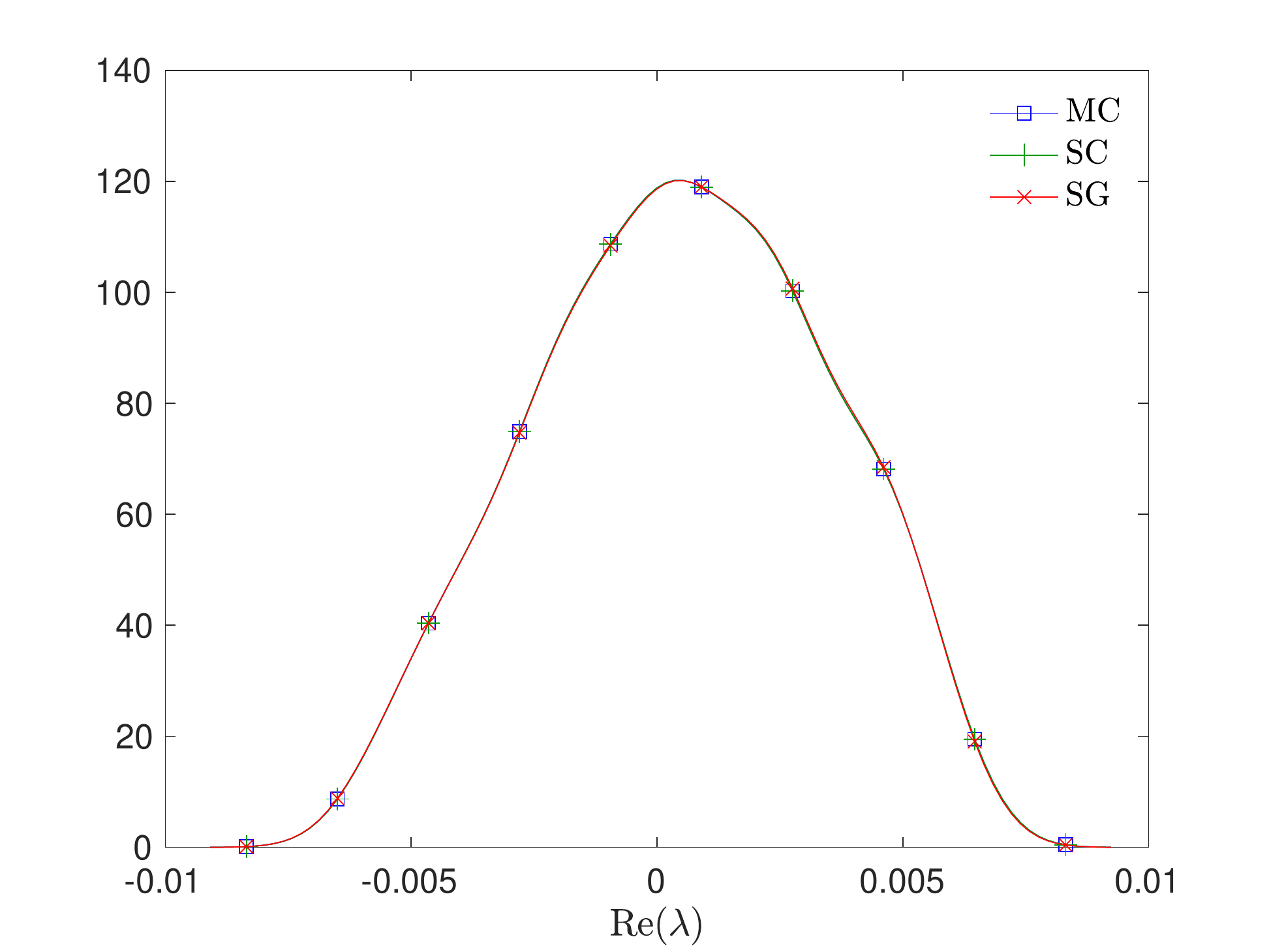}
\end{tabular}
\caption{Plots of the pdf estimate of the real part of the rightmost
eigenvalue obtained using Monte Carlo (MC), stochastic collocation (SC) and
stochastic Galerkin method (SG) for the expansion flow around a symmetric step
with $CoV=1\%$ (left) and $CoV=10\%$ (right).}%
\label{fig:symstep-pdf-Re}%
\end{figure}

For the second example, we considered an expansion flow
around a symmetric step. The domain and its discretization are shown in
Figure~\ref{fig:mesh-symstep}. The spatial discretization used a uniform grid
with $976$ \textbf{\textit{Q}}$_{2}-$\textbf{\textit{P}}$_{-1}$
finite\ elements, which provided a stable discretization for the rectangular
grid, see~\cite[p.~139]{Elman-2014-FEF}.
There were $8338$ velocity and $2928$ pressure degrees of freedom. For the
viscosity we considered a random field with affine dependence on the random
variables $\xi$ given as
\begin{equation}
\nu(x,\xi)=\nu_{1}+\sigma_{\nu}\sum_{\ell=2}^{n_{\nu}}\nu_{\ell}(x)\,\xi
_{\ell-1},\label{eq:viscosity-KL}%
\end{equation}
where $\nu_{1}$ is the mean and $\sigma_{\nu}=CoV\cdot\nu_{1}$ the standard
deviation of the viscosity, $n_{\nu}=m_{\xi}+1$, and $\nu_{\ell+1}%
=\sqrt{3\lambda_{\ell}}v_{\ell}(x)$ with $\left\{  \left(  \lambda_{\ell
},v_{\ell}(x)\right)  \right\}  _{\ell=1}^{m_{\xi}}$ are the eigenpairs of the
eigenvalue problem associated with the covariance kernel of the random field.
As in the previous example, we used the values $CoV=1\%$, and $10\%$. We
considered the same form of the covariance kernel as
in~(\ref{eq:covariance_kernel}),
\[
C\left(  X_{1},X_{2}\right)  =\exp\left(  -\frac{\left\vert x_{2}%
-x_{1}\right\vert }{L_{x}}-\frac{\left\vert y_{2}-y_{1}\right\vert }{L_{y}%
}\right)  ,
\]
and the correlation lengths were set to $12.5\%$ of the width and $25\%$\ of
the height of the domain. We assume that the random variables $\left\{
\xi_{\ell}\right\}  _{\ell=1}^{m_{\xi}}$ follow a uniform distribution over
$(-1,1)$. We note that~(\ref{eq:viscosity-KL}) can be viewed as a special case
of~(\ref{eq:viscosity}), which consists of only linear terms of$~\xi$. For the
parametrization of viscosity by~(\ref{eq:viscosity-KL}) we used the same
stochastic dimension~$m_{\xi}$ and degree of polynomial expansion$~p$ as in
the previous example: $m_{\xi}=2$\ and $p=3$, so then $n_{\xi}=10$ and
$n_{\nu}=m_{\xi}+1=3$. We used $1\times10^{3}$ samples for the Monte Carlo
method
and Smolyak sparse grid with Gauss-Legendre quadrature points and grid
level$~4$ for collocation. With these settings, the size of $\left\{  H_{\ell
}\right\}  _{\ell=1}^{n_{\nu}}$\ in~(\ref{eq:H}) was $10\times10\times3$\ with
$34$ nonzeros, and there were $n_{q}=29$ points on the sparse grid.
\textcolor{black}{As a consequence, the size of the stochastic Galerkin matrices is $112660$,
and the matrix associated with the Jacobian has in this case a block-sparse structure see, e.g.,~\cite[p.~88]{LeMaitre-2010-SMU}.}
For the solution of the Navier--Stokes problem we used the hybrid strategy
with $20$ steps of Picard iteration followed by at most $20$ steps of Newton
iteration. We used mean viscosity $\nu_{1}=4.5455\times10^{-3}$, which
corresponds to Reynolds number $Re=220$, and the rightmost eigenvalue is
$5.7963\times10^{-4}$ (the second largest eigenvalue is $-8.2273\times10^{-2}%
$), see the right panel in Figure~\ref{fig:lambdaM}.
Figure~\ref{fig:symstepMC} displays Monte Carlo realizations of the $25$
eigenvalues with the largest real part. As in the previous example, it can be
seen that the rightmost eigenvalue is relatively less sensitive to
perturbation comparing to the other eigenvalues, and it can be easily
identified in all runs of a sampling method. Figure~\ref{fig:symstep-pdf-Re}
displays the probability density function (pdf) estimates of the rightmost
eigenvalue obtained directly by Monte Carlo, the stochastic collocation and
stochastic Galerkin methods, for which the estimates were obtained using
\textsc{Matlab} function \texttt{ksdensity} for sampled gPC\ expansions. We
can see a good agreement of the plots in the left column corresponding to
$CoV=1\%$ and in the right column corresponding to $CoV=10\%$.\ In
Table~\ref{tab:symstep_gPC} we tabulate the coefficients of the gPC\ expansion
of the rightmost eigenvalue computed using the stochastic collocation and
stochastic Galerkin methods. A good agreement of coefficients can be seen, in
particular for coefficients with larger values. Finally, we examine the
inexact line-search Newton iteration from Algorithm~\ref{alg:line_search}. For
the line-search method, we used the same setup as before with $\rho=0.9$ and
$c=0.25$. The initial guess is set using the rightmost eigenvalue and
corresponding eigenvector of the eigenvalue
problem~(\ref{eq:eigenproblem-mean}), just with no imaginary part,
concatenated by zeros. The nonlinear iteration terminates when the norm of the
residual $\left\Vert \widehat{r}_{n}\right\Vert _{2}<10^{-10}$.\ The linear
systems in Line~\ref{ln:lin_sys} of Algorithm~\ref{alg:line_search} are solved
using right-preconditioned GMRES method as in the complex case.\ However,
since the eigenvalue is real, the generalized eigenvalue problem as written in
eq.~(\ref{eq:sep_param_eig2}) has the (usual) form given by
eq.~(\ref{eq:param_eig}) and all algorithms formulated in this paper can be
adapted by simply dropping the the components corresponding to the imaginary
part of the eigenvalue problem: for example, the constraints mean-based
preconditioner (Algorithm~\ref{alg:cMB}), and
specifically~(\ref{eq:M-preconditioner}) reduces to
\[
\mathcal{M}_{\text{cMB}}=\left[
\begin{array}
[c]{cc}%
K_{1}-\epsilon_{\operatorname{Re}}\mu_{\operatorname{Re}}M_{\sigma} &
-M_{\sigma}w_{\operatorname{Re}}^{(0)}\\
-w_{\operatorname{Re}}^{(0)T} & 0
\end{array}
\right]  .
\]
and in the mean-based preconditioner (Algorithm~\ref{alg:MB}) we also
modified~(\ref{eq:MB-preconditioner}) as
\[
\mathcal{M}_{\text{cMB}}=\left[
\begin{array}
[c]{cc}%
K_{1}-\epsilon_{\operatorname{Re}}\mu_{\operatorname{Re}}I & 0\\
0 & w_{\operatorname{Re}}^{(0)T}\left(  K_{1}-\epsilon_{\operatorname{Re}}%
\mu_{\operatorname{Re}}I\right)  ^{-1}w_{\operatorname{Re}}^{(0)}%
\end{array}
\right]  ,
\]
that is we used$~I$ instead of$~M_{\sigma}$ in the shift of the matrix$~K_{1}%
$. We also adapted the constraint hierarchical Gauss-Seidel preconditioner
(Algorithm~\ref{alg:chGS}--\ref{alg:chGS_cont}), which was used as before
without truncation of the matrix-vector multiplications and also with
truncation, setting $p_{t}=2$, as discussed in
Section~\ref{sec:preconditioners}. For the mean-based preconditioner we used
$\epsilon_{\operatorname{Re}}=0.97$, but the preconditioner appeared to be far
less sensitive to a specific value of$~\epsilon_{\operatorname{Re}}$, and
$\epsilon_{\operatorname{Re}}=1$ otherwise. We note that this way the
algorithms are still formulated for a generalized nonsymmetric eigenvalue
problem unlike in~\cite{Lee-2018-IMS}, where we studied symmetric problems and
in implementation we used a Cholesky factorization of the mass matrix in order
to formulate a standard eigenvalue problem. From the results in
Table~\ref{tab:gmres_symstep} we see that for all preconditioners the overall
number of nonlinear steps and GMRES\ iterations increases for larger $CoV$,
however all variants of the constraint preconditioner outperform the
mean-based preconditioner the total number of iterations remains relatively
small. Next, the performance with constraint preconditioners seem not improve
with the updating discussed in Section~\ref{sec:SMW}, which is likely since
the numbers of iterations are already low. Finally, using the constraint
hierarchical Gauss-Seidel preconditioner reduces the number of
GMRES\ iterations, which is slightly more significant for larger values of
$CoV$. \ The computational cost of preconditioner may be reduced by using the
truncation of the matrix-vector multiplications; specifically we see that the
overall iteration counts with and without the truncation are the same.

\begin{table}[b]
\caption{The $10$ coefficients of the gPC expansion of the rightmost
eigenvalue for the expansion flow around a symmetric step problem with
$CoV=1\%$ and $10\%$ computed using stochastic collocation (SC), and
stochastic Galerkin method (SG). Here $d$ is the polynomial degree and $k$ is
the index of basis function in expansion~(\ref{eq:sol_mat}).}%
\label{tab:symstep_gPC}
\begin{center}
{\footnotesize \renewcommand{\arraystretch}{1.3}
\begin{tabular}
[c]{|c|c|cc|cc|}\hline
$d$ & $k$ & SC & SG & SC & SG\\\hline
&  & \multicolumn{2}{|c|}{$CoV=1\%$} & \multicolumn{2}{|c|}{$CoV=10\%$%
}\\\hline
0 & 1 & 5.7873E-04 & 5.7873E-04 & 4.8948E-04 & 4.8927E-04\\\hline
\multirow{2}{*}{1} & 2 & -1.5948E-04 & -1.5948E-04 & -1.5890E-03 &
-1.5877E-03\\
& 3 & -2.3689E-04 & -2.3689E-04 & -2.3619E-03 & -2.3605E-03\\\hline
\multirow{3}{*}{2} & 4 & -2.4179E-07 & -2.6041E-07 & -2.4472E-05 &
-2.6501E-05\\
& 5 & -8.2562E-07 & -8.7937E-07 & -8.3136E-05 & -8.8911E-05\\
& 6 & -5.6059E-07 & -5.9203E-07 & -5.6429E-05 & -5.9831E-05\\\hline
\multirow{4}{*}{3} & 7 & 7.7918E-10 & 8.2134E-10 & 5.7810E-07 & 8.5057E-07\\
& 8 & 2.5941E-09 & 3.9327E-09 & 2.8439E-06 & 4.022E-06\\
& 9 & 3.8788E-09 & 5.5168E-09 & 4.0315E-06 & 5.6217E-06\\
& 10 & 1.3002E-09 & 2.2685E-09 & 1.6668E-06 & 2.3171E-06\\\hline
\end{tabular}
}
\end{center}
\end{table}

\begin{table}[b]
\caption{The number of GMRES iterations in each step of inexact line-search
Newton method (Algorithm~\ref{alg:line_search}) for computing the rightmost
eigenvalue and corresponding eigenvectors of the expansion flow around a
symmetric step problem with $CoV=1\%$ (left) and $10\%$ (right) and with the
stopping criteria $\Vert r_{n}\Vert_{2}<10^{-10}$ and different choices of
preconditioners: mean-based (MB) from Algorithm~\ref{alg:MB}, constraint
mean-based (cMB) from Algorithm~\ref{alg:cMB} and its updated variant
(cMB$_{u}$) from Section~\ref{sec:SMW}, and the constraint hierarchical
Gauss-Seidel preconditioner (chGS) from Algorithm~\ref{alg:chGS}%
--\ref{alg:chGS_cont} and also with truncation, setting $p_{t}=2$ (chGS$_{2}%
$). }%
\label{tab:gmres_symstep}
\begin{center}
{\footnotesize \renewcommand{\arraystretch}{1.3}
\begin{tabular}
[c]{|c|ccccc|ccccc|}\hline
& \multicolumn{5}{|c|}{$CoV=1\%$} & \multicolumn{5}{c|}{$CoV=10\%$}\\\hline
step & MB & cMB & cMB$_{u}$ & chGS & chGS$_{2}$ & MB & cMB & cMB$_{u}$ &
chGS & chGS$_{2}$\\\hline
1 & 19 & 4 & 4 & 2 & 2 & 23 & 6 & 6 & 3 & 3\\
2 & 17 & 4 & 4 & 3 & 3 & 20 & 6 & 6 & 4 & 4\\
3 & 15 & 3 & 3 & 3 & 3 & 19 & 6 & 6 & 4 & 4\\
4 &  &  &  &  &  & 15 & 5 & 5 & 4 & 4\\
5 &  &  &  &  &  & 14 & 5 & 5 & 3 & 3\\
6 &  &  &  &  &  & 23 & 8 & 8 & 5 & 5\\\hline
\end{tabular}
}
\end{center}
\end{table}

\section{Conclusion}

\label{sec:conclusion} We studied inexact stochastic Galerkin methods for
linear stability analysis of dynamical systems with parametric uncertainty and
a specific application: the Navier\textendash Stokes equation with stochastic
viscosity. The model leads to a generalized eigenvalue problem with a
symmetric mass matrix and nonsymmetric stiffness matrix, which was given by an
affine expansion obtained from a stochastic expansion of the viscosity. For
the assesment of linear stability we were interested in characterizing the
rightmost eigenvalue using the generalized polynomial chaos expansion. Since
the eigenvalue of interest may be complex, we considered separated
representation of the real and imaginary parts of the associated eigenpair.
The algorithms for solving the eigenvalue problem were formulated on the basis
of line-search Newton iteration, in which the associated linear systems were
solved using right-preconditioned GMRES\ method. We proposed several
preconditioners: the mean-based preconditioner, the constraint mean-based
preconditioner, and the constraint hierarchical Gauss-Seidel preconditioner.
For the two constraint preconditioners we also proposed an updated version,
which adapts the preconditioners to the linear system using
Sherman-Morrison-Woodburry formula \textcolor{black}{after} each step of Newton
iteration. We studied two model problems: one when the rightmost eigenvalue is
given by a complex conjugate pair and another when the eigenvalue is real. The
overall iteration count of GMRES\ with the constraint preconditioners was
smaller compared to the mean-based preconditioner, and the constraint
preconditioners were also less sensitive to the value of $CoV$. Also we found
that updating the constraint preconditioner after each step of Newton
iteration and using the off-diagonal blocks in the action of the constraint
hierarchical Gauss-Seidel preconditioner may further decrease the overall
iteration count, in particular when the rightmost eigenvalue is complex.
Finally, for both problems the probability density function estimates of the
rightmost eigenvalue closely matched the estimates obtained using the
stochastic collocation and also with the direct Monte Carlo simulation.

\paragraph{Acknowledgement}

\textcolor{black}{We would like to thank the anonymous reviewers for their insightful 
suggestions.}
Bed\v{r}ich Soused\'{\i}k would \textcolor{black}{also} like to thank
Professors Pavel Burda, Howard Elman, Vladim\'{i}r Janovsk\'{y} and Jaroslav Novotn\'{y} for many 
discussions
about bifurcation analysis and Navier\textendash Stokes equations that 
inspired 
this work.

\appendix

\section{Computations in inexact Newton iteration}

\label{sec-app:ini}The main component of a Krylov subspace method, such as
GMRES, is the computation of a matrix-vector product. First, we note that the
algorithms make use of the identity
\begin{equation}
\left(  V\otimes W\right)  \operatorname{vec}\left(  U\right)
=\operatorname{vec}\left(  WUV^{T}\right)  . \label{eq:mat-vec-equiv}%
\end{equation}
Let us write a product with Jacobian matrix from~(\ref{eq:Newton}) as
\[
\widehat{DJ}(\bar{v}_{\operatorname{Re}}^{(n)},\bar{v}_{\operatorname{Im}%
}^{(n)},\bar{\lambda}_{\operatorname{Re}}^{(n)},\bar{\lambda}%
_{\operatorname{Im}}^{(n)})%
\begin{bmatrix}
\delta\bar{v}_{\operatorname{Re}}\\
\delta\bar{v}_{\operatorname{Im}}\\
\delta\bar{\lambda}_{\operatorname{Re}}\\
\delta\bar{\lambda}_{\operatorname{Im}}%
\end{bmatrix}
,
\]
where
\begin{equation}
\widehat{DJ}(\bar{v}_{\operatorname{Re}}^{(n)},\bar{v}_{\operatorname{Im}%
}^{(n)},\bar{\lambda}_{\operatorname{Re}}^{(n)},\bar{\lambda}%
_{\operatorname{Im}}^{(n)})=\left[
\begin{array}
[c]{cccc}%
A_{\operatorname{Re}} & A_{\operatorname{Im}} & B_{\operatorname{Re}} &
B_{\operatorname{Im}}\\
-\frac{1}{2}C_{\operatorname{Re}} & -\frac{1}{2}C_{\operatorname{Im}} & 0 & 0
\end{array}
\right]  , \label{eq:jac-scheme}%
\end{equation}
with~$A_{\operatorname{Re}}$, $A_{\operatorname{Im}}$, $B_{\operatorname{Re}}%
$, $B_{\operatorname{Im}}$\ and~$C_{\operatorname{Re}}$, $C_{\operatorname{Im}%
}$ denoting the matrices in~(\ref{eq:jac_Fu_1})--(\ref{eq:jac_Gl_2}). Then,
using~(\ref{eq:U}) and~(\ref{eq:mat-vec-equiv}), the matrix-vector product
entails evaluating%
\begin{align}
A_{\operatorname{Re}}\delta\bar{v}_{\operatorname{Re}}  &  =%
\begin{bmatrix}
\mathbb{E}\left[  \Psi\Psi^{T}\!\otimes\!K\right]  \!-\!\mathbb{E}%
[(\bar{\lambda}_{\operatorname{Re}}^{(n)T}\Psi)\Psi\Psi^{T}\!\otimes
\!M_{\sigma}]\\
\vspace{-2.5mm}\\
-\mathbb{E}[(\bar{\lambda}_{\operatorname{Im}}^{(n)T}{}\Psi)\Psi\Psi
^{T}\!\otimes\!M_{\sigma}]
\end{bmatrix}
\delta\bar{v}_{\operatorname{Re}}=\label{eq:jac_A_Re}\\
&  =%
\begin{bmatrix}
\operatorname{vec}\left(  \sum_{\ell=1}^{n_{\nu}}K_{\ell}\delta\bar
{V}_{\operatorname{Re}}H_{\ell}^{T}\right)  \!-\!\operatorname{vec}\left(
\sum_{\ell=1}^{n_{\xi}}\lambda_{\operatorname{Re},\ell}^{(n)}M_{\sigma}%
\delta\bar{V}_{\operatorname{Re}}H_{\ell}^{T}\right) \\
\vspace{-2.5mm}\\
-\operatorname{vec}\left(  \sum_{\ell=1}^{n_{\xi}}\lambda_{\operatorname{Im}%
,\ell}^{(n)}M_{\sigma}\delta\bar{V}_{\operatorname{Re}}H_{\ell}^{T}\right)
\end{bmatrix}
,\nonumber\\
A_{\operatorname{Im}}\delta\bar{v}_{\operatorname{Im}}  &  =%
\begin{bmatrix}
\mathbb{E}[(\bar{\lambda}_{\operatorname{Im}}^{(n)T}{}\Psi)\Psi\Psi
^{T}\!\otimes\!M_{\sigma}]\\
\vspace{-2.5mm}\\
\mathbb{E}[\Psi\Psi^{T}\!\otimes\!K]\!-\!\mathbb{E}[(\bar{\lambda
}_{\operatorname{Re}}^{(n)T}{}\Psi)\Psi\Psi^{T}\!\otimes\!M_{\sigma}]
\end{bmatrix}
\delta\bar{v}_{\operatorname{Im}}=\label{eq:jac_A_Im}\\
&  =%
\begin{bmatrix}
\operatorname{vec}\left(  \sum_{\ell=1}^{n_{\xi}}\lambda_{\operatorname{Im}%
,\ell}^{(n)}M_{\sigma}\delta\bar{V}_{\operatorname{Re}}H_{\ell}^{T}\right) \\
\vspace{-2.5mm}\\
\operatorname{vec}\left(  \sum_{\ell=1}^{n_{\nu}}K_{\ell}\delta\bar
{V}_{\operatorname{Re}}H_{\ell}^{T}\right)  \!-\!\operatorname{vec}\left(
\sum_{\ell=1}^{n_{\xi}}\lambda_{\operatorname{Re},\ell}^{(n)}M_{\sigma}%
\delta\bar{V}_{\operatorname{Re}}H_{\ell}^{T}\right)
\end{bmatrix}
,\nonumber
\end{align}%
\begin{align}
B_{\operatorname{Re}}\delta\bar{\lambda}_{\operatorname{Re}}  &  =%
\begin{bmatrix}
-\mathbb{E}[\Psi^{T}\otimes(\Psi\Psi^{T}\!\otimes\!M_{\sigma})]\bar
{v}_{\operatorname{Re}}^{(n)}\\
\vspace{-2.5mm}\\
-\mathbb{E}[\Psi^{T}\otimes(\Psi\Psi^{T}\!\otimes\!M_{\sigma})]\bar
{v}_{\operatorname{Im}}^{(n)}%
\end{bmatrix}
\delta\bar{\lambda}_{\operatorname{Re}}=%
\begin{bmatrix}
-\operatorname{vec}\left(  \sum_{\ell=1}^{n_{\xi}}\delta\lambda
_{\operatorname{Re},\ell}M_{\sigma}\bar{V}_{\operatorname{Re}}^{(n)}H_{\ell
}^{T}\right) \\
\vspace{-2.5mm}\\
-\operatorname{vec}\left(  \sum_{\ell=1}^{n_{\xi}}\delta\lambda
_{\operatorname{Re},\ell}M_{\sigma}\bar{V}_{\operatorname{Im}}^{(n)}H_{\ell
}^{T}\right)
\end{bmatrix}
,\label{eq:jac_B_Re}\\
B_{\operatorname{Im}}\bar{\lambda}_{\operatorname{Im}}  &  =%
\begin{bmatrix}
\mathbb{E}[\Psi^{T}\otimes(\Psi\Psi^{T}\!\otimes\!M_{\sigma})]\bar
{v}_{\operatorname{Im}}^{(n)}\\
\vspace{-2.5mm}\\
-\mathbb{E}[\Psi^{T}\otimes(\Psi\Psi^{T}\!\otimes\!M_{\sigma})]\bar
{v}_{\operatorname{Re}}^{(n)}%
\end{bmatrix}
\delta\bar{\lambda}_{\operatorname{Im}}=%
\begin{bmatrix}
\operatorname{vec}\left(  \sum_{\ell=1}^{n_{\xi}}\delta\lambda
_{\operatorname{Im},\ell}M_{\sigma}\bar{V}_{\operatorname{Im}}^{(n)}H_{\ell
}^{T}\right) \\
\vspace{-2.5mm}\\
-\operatorname{vec}\left(  \sum_{\ell=1}^{n_{\xi}}\delta\lambda
_{\operatorname{Im},\ell}M_{\sigma}\bar{V}_{\operatorname{Re}}^{(n)}H_{\ell
}^{T}\right)
\end{bmatrix}
, \label{eq:jac_B_Im}%
\end{align}%
\begin{align}
-\frac{1}{2}C_{\operatorname{Re}}\delta\bar{v}_{\operatorname{Re}}  &  =-%
\begin{bmatrix}
\mathbb{E}[\Psi\otimes(\bar{v}_{\operatorname{Re}}^{(n)T}\Psi\Psi^{T}%
\!\otimes\!I_{n_{x}})]\\
\vspace{-2.5mm}\\
0
\end{bmatrix}
\delta\bar{v}_{\operatorname{Re}}=-%
\begin{bmatrix}%
\begin{bmatrix}
\bar{v}_{\operatorname{Re}}^{(n)T}(H_{1}\otimes I_{n_{x}})\delta\bar
{v}_{\operatorname{Re}}\\
\vdots\\
\bar{v}_{\operatorname{Re}}^{(n)T}(H_{n_{\xi}}\otimes I_{n_{x}})\delta\bar
{v}_{\operatorname{Re}}%
\end{bmatrix}
\\
\vspace{-2.5mm}\\
0
\end{bmatrix}
,\label{eq:jac_C_Re}\\
-\frac{1}{2}C_{\operatorname{Im}}\delta\bar{v}_{\operatorname{Im}}  &  =-%
\begin{bmatrix}
0\\
\vspace{-2.5mm}\\
\mathbb{E}[\Psi\otimes(\bar{v}_{\operatorname{Im}}^{(n)T}\Psi\Psi^{T}%
\!\otimes\!I_{n_{x}})]
\end{bmatrix}
\delta\bar{v}_{\operatorname{Im}}=-%
\begin{bmatrix}
0\\
\vspace{-2.5mm}\\%
\begin{bmatrix}
\bar{v}_{\operatorname{Im}}^{(n)T}(H_{1}\otimes I_{n_{x}})\delta\bar
{v}_{\operatorname{Im}}\\
\vdots\\
\bar{v}_{\operatorname{Im}}^{(n)T}(H_{n_{\xi}}\otimes I_{n_{x}})\delta\bar
{v}_{\operatorname{Im}}%
\end{bmatrix}
\end{bmatrix}
, \label{eq:jac_C_Im}%
\end{align}
and the right-hand side of~(\ref{eq:Newton}) is evaluated using
\begin{align*}
F^{(n)}  &  =%
\begin{bmatrix}
\mathbb{E}[\Psi\Psi^{T}\!\otimes\!K]\bar{v}_{\operatorname{Re}}^{(n)}%
\!-\!\mathbb{E}[(\bar{\lambda}_{\operatorname{Re}}^{(n)T}{}\Psi)\Psi\Psi
^{T}\!\otimes\!M_{\sigma}]\bar{v}_{\operatorname{Re}}^{(n)}\!+\!\mathbb{E}%
[(\bar{\lambda}_{\operatorname{Im}}^{(n)T}{}\Psi)\Psi\Psi^{T}\!\otimes
\!M_{\sigma}]\bar{v}_{\operatorname{Im}}^{(n)}\\
\vspace{-2.5mm}\\
\mathbb{E}[\Psi\Psi^{T}\!\otimes\!K]\bar{v}_{\operatorname{Im}}^{(n)}%
\!-\!\mathbb{E}[(\bar{\lambda}_{\operatorname{Im}}^{(n)T}{}\Psi)\Psi\Psi
^{T}\!\otimes\!M_{\sigma}]\bar{v}_{\operatorname{Re}}^{(n)}\!-\!\mathbb{E}%
[(\bar{\lambda}_{\operatorname{Re}}^{(n)T}{}\Psi)\Psi\Psi^{T}\!\otimes
\!M_{\sigma}]\bar{v}_{\operatorname{Im}}^{(n)}%
\end{bmatrix}
=\\
&  =%
\begin{bmatrix}
\operatorname{vec}\left(  \sum_{\ell=1}^{n_{\nu}}K_{\ell}\delta\bar
{V}_{\operatorname{Re}}^{(n)}H_{\ell}^{T}\right)  \!-\!\operatorname{vec}%
\left(  \sum_{\ell=1}^{n_{\xi}}\lambda_{\operatorname{Re},\ell}^{(n)}%
M_{\sigma}\bar{V}_{\operatorname{Re}}^{(n)}H_{\ell}^{T}\right)
\!+\!\operatorname{vec}\left(  \sum_{\ell=1}^{n_{\xi}}\lambda
_{\operatorname{Im},\ell}^{(n)}M_{\sigma}\bar{V}_{\operatorname{Im}}%
^{(n)}H_{\ell}^{T}\right) \\
\vspace{-2.5mm}\\
\operatorname{vec}\left(  \sum_{\ell=1}^{n_{\nu}}K_{\ell}\delta\bar
{V}_{\operatorname{Im}}^{(n)}H_{\ell}^{T}\right)  \!-\!\operatorname{vec}%
\left(  \sum_{\ell=1}^{n_{\xi}}\lambda_{\operatorname{Im},\ell}^{(n)}%
M_{\sigma}\bar{V}_{\operatorname{Re}}^{(n)}H_{\ell}^{T}\right)
\!-\!\operatorname{vec}\left(  \sum_{\ell=1}^{n_{\xi}}\lambda
_{\operatorname{Re},\ell}^{(n)}M_{\sigma}\bar{V}_{\operatorname{Im}}%
^{(n)}H_{\ell}^{T}\right)
\end{bmatrix}
,
\end{align*}
and
\[
G^{(n)}=%
\begin{bmatrix}
\mathbb{E}[\Psi\otimes((\bar{v}_{\operatorname{Re}}^{(n)T}(\Psi\Psi^{T}\otimes
I_{n_{x}})\bar{v}_{\operatorname{Re}}^{(n)})\!-\!1)]\\
\vspace{-2.5mm}\\
\mathbb{E}[\Psi\otimes((\bar{v}_{\operatorname{Im}}^{(n)T}(\Psi\Psi^{T}\otimes
I_{n_{x}})\bar{v}_{\operatorname{Im}}^{(n)})\!-\!1)]
\end{bmatrix}
,
\]
where, using $\star$ for either $\operatorname{Re}$ or $\operatorname{Im}$,
the $i$th row of$~G^{(n)}$ is
\begin{align*}
\left[  G^{(n)}\right]  _{i}  &  =\mathbb{E}[\psi_{i}(\bar{v}_{\star}^{(n)}%
{}^{T}(\Psi\Psi^{T}\otimes I_{n_{x}})\bar{v}_{\star}^{(n)})-\psi_{i}],\\
&  =\bar{v}_{\star}^{(n)}{}^{T}\mathbb{E}[\psi_{i}\Psi\Psi^{T}\otimes
I_{n_{x}}]\bar{v}_{\star}^{(n)}-\delta_{1i},
\end{align*}
and the first term above is evaluated as
\[
\bar{v}_{\star}^{(n)}{}^{T}\mathbb{E}[\psi_{i}\Psi\Psi^{T}\otimes I_{n_{x}%
}]\bar{v}_{\star}^{(n)}=\bar{v}_{\star}^{(n)}{}^{T}(H_{i}\otimes I_{n_{x}%
})\bar{v}_{\star}^{(n)},
\]
or, denoting the trace operator by$~\operatorname{tr}$, this term can be also
evaluated as
\[
\bar{v}_{\star}^{(n)}{}^{T}\mathbb{E}[\psi_{i}\Psi\Psi^{T}\otimes I_{n_{x}%
}]\bar{v}_{\star}^{(n)}=\operatorname{tr}(\bar{V}_{\star}^{(n)}H_{i}\bar
{V}_{\star}^{(n)}{}^{T})=\operatorname{tr}(\bar{V}_{\star}^{(n)}{}^{T}\bar
{V}_{\star}^{(n)}H_{i}).
\]

{\color{black}
\section{Parameters used in numerical experiments}

In addition to the description in Section~\ref{sec:numerical}, we provide in Table~\ref{tab:parameters} an
overview of the main parameters used in the numerical experiments. 
Besides that, we used the following settings in both experiments. The gPC parameters:
$m_\xi=2$, $p=3$, $n_\xi=10$; for the sampling methods: $n_{MC}=1\times 10^3$, $n_q=29$;
for the inexact Newton iteration: $\rho=0.9$, $c=0.25$, stopping criterion $\left\Vert \widehat{r}_{n}\right\Vert _{2}<10^{-10}$;
for the preconditioners: $\epsilon_{\operatorname{Re}}=\epsilon_{\operatorname{Im}}=0.97$ (the mean-based preconditioner) 
and $\epsilon_{\operatorname{Re}}=\epsilon_{\operatorname{Im}}=1$ (otherwise). 

\begin{table}[h]
\caption{The main parameters used in the numerical experiments.}%
\label{tab:parameters}
\begin{center}
{\footnotesize \renewcommand{\arraystretch}{1.3}
\begin{tabular}
[c]{|c|c|c|}\hline
 & Section~\ref{sec:obstacle} & Section~\ref{sec:expansion} \\ \hline
 problem & Flow around an obstacle & Expansion flow around a symmetric step \\ \hline
FEM & \textbf{\textit{Q}}$_{2}-$\textbf{\textit{Q}}$_{1}$ & \textbf{\textit{Q}}$_{2}-$\textbf{\textit{P}}$_{-1}$  \\
nelem/$n_u$/$n_p$ & 1008/8416/1096 & 976/8338/2928 \\ 
$Re$ & 373 & 220 \\ 
$\lambda$ & $0.0085\pm2.2551i$ & $5.7963\times 10^{-4}$ \\ 
$n_\nu$ & 28 & 3 \\
quadrature (in SC) & Gauss-Hermite & Gauss-Legendre \\ \hline
\multicolumn{3}{|c|}{Solving the Navier\textendash Stokes problem (see~\cite{Sousedik-2016-SGM} for details):} \\ \hline
max Picard steps & 6 & 20 \\
max Newton steps & 15 & 20 \\
nltol & $10^{-8}$ & $10^{-8}$ \\ \hline
\end{tabular}
}
\end{center}
\end{table}
}

\clearpage
\bibliographystyle{plain}
\bibliography{eig_ns.bib}

\end{document}